\documentclass[english, 12pt]{amsart}

\usepackage{scrtime}

\usepackage{tikz}
\usetikzlibrary{cd}
\usepackage{aliascnt}
\usepackage{mathrsfs}
\usepackage[all,poly,knot]{xy}
\usepackage{hyperref}
\usepackage{comment}  
\usepackage{csquotes}   
\usepackage{amssymb,amsbsy,amsmath,amsfonts,amssymb,amscd,
	graphics,color,footmisc,fancyhdr,multicol,fancybox,
	graphicx,mathrsfs,rotating,ifthen,wasysym}

\usepackage{lipsum}

\usepackage{times}
\usepackage{mathtools}
\usepackage{pdfpages}
\usepackage{multirow}
\usepackage{nccrules}
\usepackage{textcomp}
\usepackage{comment}
\usepackage{framed}
\usepackage[all]{xy}
\usepackage{graphicx}
\usepackage{pgf,tikz}
\usepackage{mathrsfs}
\usetikzlibrary{arrows}
\usepackage{lipsum}
\usepackage{mathtools}

\DeclareMathOperator{\dif}{\text{\normalfont d}}

\def\log{\mathrm{log}\,}
\newcommand{\norm}[1]{\lVert#1\rVert}

\theoremstyle{plain}

\newtheorem{thm}{Theorem}[section]

\newtheorem{lem}[thm]{{Lemma}}

\newtheorem{pro}[thm]{Proposition}

\theoremstyle{definition}

\newtheorem{defi}[thm]{{Definition}}

\newtheorem{ques}[thm]{{Question}}
\newtheorem{obs}[thm]{{Observation}} 

\theoremstyle{remark}
\newtheorem{rmk}[thm]{Remark}

\numberwithin{equation}{section}

\usepackage[top=2cm, bottom=2cm, left=2.5cm, right=2.5cm]{geometry}

\usepackage[hyperpageref]{backref}

\usepackage{xcolor}
\hypersetup{
colorlinks,
linkcolor={red!50!black},
citecolor={blue!62!black},
urlcolor={blue!80!black}
}
\theoremstyle{plain}
\newcommand{\thistheoremname}{}
\newtheorem*{genericthm*}{\thistheoremname}
\newenvironment{namedthm*}[1]{\renewcommand{\thistheoremname}{#1}%
\begin{genericthm*}}
{\end{genericthm*}}

\newtheoremstyle{named}{}{}{\itshape}{}{\bfseries}{.}{.5em}{\thmnote{#3's }#1}
\theoremstyle{named}

\makeatletter
\newcommand\thankssymb[1]{\textsuperscript{\@fnsymbol{#1}}}
\makeatother

\begin{document}

	\title{
		Universal Holomorphic Maps with Slow Growth
		\\
		II. Functional Analysis Methods
}


\author{Bin Guo}

\address{Academy of Mathematics and Systems Sciences, Chinese Academy of Sciences, Beijing 100190, China}
\email{guobin181@mails.ucas.ac.cn}

\author{Song-Yan Xie }
\address{Academy of Mathematics and System Science \& Hua Loo-Keng Key Laboratory
		of Mathematics, Chinese Academy of Sciences, Beijing 100190, China; 
		School of Mathematical Sciences, University of Chinese Academy of Sciences, Beijing
100049, China}
\email{xiesongyan@amss.ac.cn}

\date{\today}


\begin{abstract}
By means of hypercyclic operator theory, 
we complement
  our previous results  
 on hypercyclic  holomorphic maps between complex Euclidean spaces having slow growth rates,
 by showing {\it abstract abundance} rather than 
 {\it explicit existence}. 
 Next, 
 we establish that,
 in the space of holomorphic maps from $\mathbb{C}^n$  to any connected Oka manifold $Y$, equipped with the compact-open topology, there exists a {\em dense} 
 subset 
 consisting of common {\em frequently hypercyclic}
 elements
 for all nontrivial translation operators. This is new even for $n=1$ and $Y=\mathbb{C}$.
\end{abstract}

\maketitle

\section{\bf Introduction}
On the  space of entire holomorphic functions $\mathcal{H}(\mathbb{C})$, equipped with the compact-open topology, the translation operator $\mathsf{T}_a\,f(\bullet) = f(\bullet+a)$ acts linearly for every  $a\in\mathbb{C}$. 
Birkhoff~\cite{Birkhoff1929} discovered that   $\mathcal{H}(\mathbb{C})$
contains {\em universal} elements $f$, in the sense that a certain countable set of translations $\{\mathsf{T}_{a_i}f\}_{i=1}^{\infty}$ is dense in $\mathcal{H}(\mathbb{C})$. In fact, for any  $\mathsf{T}_a$ with $a\neq 0$, Birkhoff can make $f\in \mathcal{H}(\mathbb{C})$ a {\em hypercyclic} element for $\mathsf{T}_a$, {\em i.e.},
the iterated orbit $\{\mathsf{T}_{a}^{i} f=\mathsf{T}_{i\cdot a}\,f\}_{i\geqslant 1}$ is dense.
This surprising phenomenon 
opened the door for a new branch of functional analysis / dynamical systems called linear chaos~\cite{MR1685272, book-bayart-matheron, book-linear-chaos}, according to which infinite dimensional linear dynamics  may exhibit  exotic patterns.

By complex analysis, a universal entire function $f$ must be transcendental, and hence its asymptotic growth rate has to be larger than that of any  polynomial. One might expect that the minimal possible growth rate of $f$  shall be  considerably larger than this obvious bound. However, a surprising result~\cite{DRuis-paper-1} shows the contrary.

\begin{thm}[Du\u{\i}os Ruis]
	\label{thm-Duios-Ruis}
	For a continuous increasing function
	$\phi: \mathbb{R}_{\geqslant 0} \rightarrow \mathbb{R}_{> 0}$
	growing faster than any polynomials
	\begin{equation}\label{growth phi}
	\lim_{r\rightarrow \infty}
		\phi(r)/r^N
		=\infty
		\qquad
		{\scriptstyle
			(\forall\, N\, \geqslant\, 1),
		}
	\end{equation}
	 there exists a universal entire function $f$ with slow growth 
	\[
	|f(z)|\leqslant
		\phi(|z|)
		\qquad
		{\scriptstyle
			(\forall\,z\,\in\, \mathbb{C})
		}.   
	\]
	\qed
\end{thm}

The original exposition in~\cite{DRuis-paper-1} is not easy to decipher. A breakthrough~\cite{Chan-Shapiro-91} of Chan and Shapiro contributed an alternative proof by functional analysis, and their method has been gradually absorbed, simplified,  and  
became standard in  hypercyclic operator theory~(cf.~\cite{desch1997hypercyclic} and the references therein).

In a similar flavor,
 Dinh and Sibony~\cite[Problem 9.1]{Dinh-Sibony-list} raised a problem of seeking
 minimal growth rate of universal  meromorphic functions $g: \mathbb{C} \rightarrow \mathbb{CP}^1$ in terms of the Nevanlinna-Shimizu-Ahlfors  characteristic function
$$T_g(r)=\int_1^r\frac{\dif t}{t}\int_{\mathbb{D}_t}g^*\omega_{\mathsf{FS}},$$
where 
$\omega_{\mathsf{FS}}$ is the Fubini-Study form on $\mathbb{CP}^1$ giving unit area $\int_{\mathbb{CP}^1}\omega_{\mathsf{FS}}=1$,
and $\mathbb{D}_t\subset \mathbb{C}$ is the disc centered at the origin with radius $t$. 
An optimal answer was obtained as

\begin{thm}[\cite{Xie-Chen-Tuan}]
	\label{chen-hunh-xie thm}
	For any positive continuous nondecreasing function 
	\begin{equation*}\label{growth varphi}
	    \psi: 
	\mathbb{R}_{\geqslant 1}
	\rightarrow
	\mathbb{R}_+
	\qquad
	\text{satisfying}
	\qquad 
	\lim_{r\rightarrow +\infty}\psi(r)
    =
    +\infty,
	\end{equation*}
	there exists some
	universal meromorphic function $g: \mathbb{C}\rightarrow \mathbb{CP}^n$ with  slow growth
	\begin{equation*}
		\label{admissible-slow-growth}
		T_g(r)
		\leqslant
		\psi(r)
		\cdot
		\log r
		\qquad
		{\scriptstyle
			(\forall\, r\geqslant\, 1).
		}
	\end{equation*}
	\qed
\end{thm}

By a constructive algorithm~\cite{Guo-Xie-1}, 
Theorems \ref{thm-Duios-Ruis}, ~\ref{chen-hunh-xie thm} have been generalized to the higher dimensional case as follows.

\begin{defi}
\label{hypercyclic}
	Let $X$ be a topological space.
	A continuous operator $\mathsf{T}:X\rightarrow X$ is called {\sl hypercyclic} if  some element $x\in X$ has dense   orbit  
	\[
	\mathsf{orb}(x,\mathsf{T}):=\{x, \mathsf{T}^{1}(x), \mathsf{T}^{2}(x), \mathsf{T}^{3}(x), \dots\}
	\]  
	in $X$. Such an $x$ is called a hypercyclic element for $\mathsf{T}$. The set that consists of all hypercyclic elements for $\mathsf{T}$
	is denoted by $\mathsf{HC}(\mathsf{T})$. 
\end{defi}

\begin{thm}[\cite{Guo-Xie-1}]\label{main in math ann}
    Given any transcendental  growth function $\phi$ growing faster than any polynomial, 
	for
	countably many directions $\{\theta_i\}_{i\geqslant 1}$ in the unit sphere $\mathbb{S}^{2n-1}\subset
	\mathbb{C}^n$, there exist
	universal  holomorphic maps
	$F: \mathbb{C}^n \rightarrow \mathbb{C}^m$ with slow growth 
	\begin{equation*}
		\label{desired growth in several variables}
		\norm{F(z)}\leqslant
		\phi(\norm{z})
		\qquad
		\qquad
		{\scriptstyle
			(\forall\,z\,\in\, \mathbb{C}^n),
		}
	\end{equation*}
	such that $F$ are hypercyclic for all  $\mathsf{T}_a$, where  $a$ are in some ray $\mathbb{R}_{+}\cdot \theta_i$. \qed
\end{thm}

\begin{thm}[\cite{Guo-Xie-1}]
    Let $n\leqslant m$ be two positive integers. 
	Let $\psi: \mathbb{R}_{\geqslant 1}\rightarrow \mathbb{R}_+$
	be a continuous nondecreasing function tending to infinity. 
    For $F\in\mathsf{Hol}(\mathbb{C}^n,\mathbb{CP}^m)$,  the associated Nevanlinna characteristic function is defined by
	\[
T_F(r)
:=
\int_{1}^r
\frac{\dif t}{t^{2n-1}}
\int_{\{\norm{z}<t\}}
F^*\omega_{\mathsf{FS}}\wedge \alpha^{n-1}
\qquad
		{\scriptstyle
			(\forall\,r\,\geqslant\, 1),}
\]
where $\omega_{\mathsf{FS}}$ is the Fubini-Study metric on $\mathbb{CP}^m$ and $\alpha:=\dif\dif^c \norm{z}^2 $ is a  positive  $(1, 1)$ form.
	Then
	for countably many
	given directions  $\{\theta_i\}_{i\geqslant 1}$ in the unit sphere $\mathbb{S}^{2n-1}\subset
	\mathbb{C}^n$, there exist
	universal entire holomorphic maps
	$f: \mathbb{C}^n \rightarrow \mathbb{CP}^m$ with slow growth 
	\begin{equation*}
	\label{admissible-slow-growth cpn}
		T_f(r)
		\leqslant
		\psi(r)
		\cdot
		\log r
		\qquad
		{\scriptstyle
			(\forall\, r\,\geqslant\, 1),
		}
	\end{equation*}
	such that  $f$ are hypercyclic for all nontrivial translations 
	$\mathsf{T}_a$,
	where  $a$ are in some ray
	$\mathbb{R}_+\cdot \theta_i$. 
	\qed
\end{thm}

The constructive approach of~\cite{Guo-Xie-1}
has the advantage of  writing explicitly universal holomorphic maps $F$ having slow growth,
from which delicate information can be read off. For instance, we can show that

\begin{thm}[\cite{Guo-Xie-1}]
	Given any transcendental growth function  $\phi$ of the shape~\eqref{growth phi},
	there exists some
	universal entire holomorphic function
	$F\in \mathcal{H}(\mathbb{C})$ with slow growth
	\[
	|F(z)|\leqslant
	\phi(|z|)
	\qquad
	\qquad
	{\scriptstyle
		(\forall\,z\,\in\, \mathbb{C}),
	}
	\]
	such that  the set $I\subset \mathbb{S}^1$ of hypercyclic directions of $F$ is uncountable.
	\qed
\end{thm}


However, our constructive method fails to tell the ``density'' of such $F$ in Theorem~\ref{main in math ann}.   
Our original motivation in this work was to understand Theorem~\ref{main in math ann} by means of hypercyclic operator theory~(cf.~\cite{book-bayart-matheron, book-linear-chaos}). Indeed, by implementing functional analysis methods (cf.~\cite{Chan-Shapiro-91, desch1997hypercyclic}), we can show that such $F$ are reasonably abundant.

\medskip

Let $\phi$ be given as \eqref{growth phi}.
 Consider the set
\[
\mathsf{S}_{\phi}\,
:=\,
\{f \in\mathsf{Hol}(\mathbb{C}^n,\mathbb{C}^m) \,:\, \norm{f(z)}\leqslant
		\phi(\norm{z})
		\,,
		{
			\forall\, z\in \mathbb{C}^n
		}  \}
\]
of entire functions having slow growth with respect to $\phi$. 
Hence
we can rephrase Theorem~\ref{main in math ann} as
$$
\mathsf{A}_\phi:=
\cap_{k\geqslant 1}\cap_{ r\in\mathbb{R}_+}\mathsf{HC}(\mathsf{T}_{r\cdot e^{[k]}})\cap \mathsf{S}_{\phi}\neq\emptyset
		\qquad
		{\scriptstyle
			(\text{given }\,e^{[k]}\,\in\, \mathbb{S}^{2n-1}\,\text{for}\,k\,\geqslant\,1).
		}
$$

Can we say more about $\mathsf{A}_\phi$? First, it is not everywhere dense, since we can take a  holomorphic map $G\in\mathsf{Hol}(\mathbb{C}^n,\mathbb{C}^m)$ with $\norm{G(1)}>\phi(1)+2$. Hence the open neighborhood 
\[
\{h\in\mathsf{Hol}(\mathbb{C}^n,\mathbb{C}^m):\sup_{z\in\overline{\mathbb{B}}(0,2)} \norm{h(z)-G(z)}<1\}
\]
of $G$
has empty intersection with $\mathsf{S}_{\phi}$.
Nevertheless, $\mathsf{A}_\phi$ is reasonably dense near $\mathbf{0}\in\mathsf{Hol}(\mathbb{C}^n,\mathbb{C}^m)$.

\medskip\noindent {\bf Theorem A.}
{\it
 For arbitrary countable vectors $\{e^{[k]}\}_{k\geqslant 1}$ in $\mathbb{S}^{2n-1}$ and for $\phi$ as~\eqref{growth phi}, 
 there exists a Banach space $\mathbf{B}$ and a continuous linear injective map $\iota: \mathbf{B}\rightarrow \mathsf{Hol}(\mathbb{C}^n,\mathbb{C}^m)$ with dense image, 
 such that
 for some $G_\delta$-dense subset $G$ of $\mathbf{B}$ there holds
 $$
 \iota(G\cap \mathbb{B}(0,1))\,
 \subset\,
 \cap_{k\geqslant 1}\cap_{ r\in\mathbb{R}_+}\mathsf{HC}(\mathsf{T}_{r\cdot e^{[k]}})\cap \mathsf{S}_{\phi},
 $$
 where 
 $\mathbb{B}(0,1)\subset \mathbf{B}$ is the unit ball centered at $0$. 
}

\medskip

Moreover, we can generalize  Theorem~A to an infinite dimensional case. The insight comes from the tradition in several complex variables: Oka's Jok\^u-Iko Principle (cf.~\cite{oka1936fonctions, MR3526579}).
We can imagine that the  countably many directions $\{e^{[k]}\}_{k\geqslant 1}$ in $\mathbb{C}^n$ are just the ``shadows'' or projections
of the infinite coordinate directions of $\mathbb{C}^{\infty}$ onto the finite dimensional subspace $\mathbb{C}^n\hookrightarrow$ $\mathbb{C}^{\infty}$, 
thus we can guess that a
certain analogous statement holds true
for ``$\mathsf{Hol}(\mathbb{C}^{\infty}, \mathbb{C}^{\infty})$''.
However, a proper definition of ``$\mathsf{Hol}(\mathbb{C}^{\infty}, \mathbb{C}^{\infty})$''
requires much effort, see Sect.~\ref{Infinite dimension analytic function}.

\smallskip\noindent{\bf Convention.} We abuse the same notation
$\mathsf{S}_{\phi}$ to denote subsets of
$\mathsf{Hol}(\mathbb{C}^n,\mathbb{C}^m)$, or
$\mathcal A^{\mathsf{w}}_{\mathsf{fin}}(\ell^2(\mathbb{C}),\ell^2(\mathbb{C}))$ (see~\eqref{analytic maps between infinitely complex spaces}),
which consist of elements $f$ satisfying
$\norm{f(z)}\leqslant
		\phi(\norm{z})$ for all possible $z$ according to the context.
Let $\mathbb{C}_c^{\mathbb{N}}$ be the subset of $ \ell^2(\mathbb{C})$ consisting of elements$(x_i)_{i\geqslant 1}$ such that all but finitely many $x_i=0$.

\medskip\noindent {\bf Theorem B.}
{\it
 For arbitrary countablely many  nonzero vectors $\{a^{[k]}\}_{k\geqslant 1}$ in $\mathbb{C}_c^{\mathbb{N}}$, the set  of common hypercyclic elements  in $\mathcal A^{\mathsf{w}}_{\mathsf{fin}}(\ell^2(\mathbb{C}),\ell^2(\mathbb{C}))$ for the translation operators
 $\{\mathsf{T}_{a^{[k]}}\}_{k\geqslant 1}$, 
 having slow growth with respect to $\phi$,  
	is nonempty
	\[
	\cap_{k\geqslant 1}\,\mathsf{HC}(\mathsf{T}_{a^{[k]}})\cap \mathsf{S}_{\phi}
	\neq\emptyset.
	\]
	Moreover, 
this  set contains certain image $\iota(\mathbb{B}(0,1)\cap G)$, where
 $\iota: \mathbf{B}\hookrightarrow \mathcal{A}^{\mathsf{w}}_{\mathsf{fin}}(\ell^2(\mathbb{C}),\ell^2(\mathbb{C}))$ is a continuous linear injective map from some Banach space $\mathbf{B}$, and where 
 $\mathbb{B}(0,1)\subset\mathbf{B}$ is the unit ball centered at $0$, and where
 $G\subset \mathbf{B}$ is some 
$G_\delta$-dense subset. 
}

\medskip
In hypercyclic operator theory, it is known that $
\cap_{a\in\mathbb{C}\setminus\{0\}}\mathsf{HC}(\mathsf{T}_a)\neq\emptyset$~(cf.~\cite[Example 11.21]{book-linear-chaos}, see also the discussion~\cite[ p.~330]{book-linear-chaos}).
A breakthrough of  Costakis and Sambarino shows moreover  the abundance.

\begin{thm}[\cite{amazing-theorem}]\label{all trans}
The set $\cap_{a\in\mathbb{C}\setminus\{0\}}\mathsf{HC}(\mathsf{T}_a)$ of common hypercyclic elements for all nontrivial translation operators    contains a $G_\delta$-dense subset of $\mathcal{H}(\mathbb{C})$. 
\qed
\end{thm}

Hence it is natural to ask whether or not 
$
\cap_{a\in\mathbb{C}\setminus\{0\}}\mathsf{HC}(\mathsf{T}_a)\cap\mathsf{S}_\phi\neq\emptyset
$
for any $\phi$ satisfying~\eqref{growth phi}. The answer is
{\em No}~\cite[Theorem~D]{Guo-Xie-1}, by Nevanlinna theory (cf.~\cite{MR3156076, MR4265173}). 

\begin{thm}[\cite{Guo-Xie-1}]\label{zero measure 1 dim}
	There exists some
	transcendental growth function   $\phi$ of the shape~\eqref{growth phi},
	such that for any entire holomorphic function
	$F\in \mathcal{H}(\mathbb{C})$ with slow growth
	\begin{equation*}
	|F(z)|\leqslant
		\phi(|z|)
		\qquad
		\qquad
		{\scriptstyle
			(\forall\,z\,\in\, \mathbb{C}),
		}
	\end{equation*}
	the set $\{a\in  \mathbb{S}^1: F \text{\,is hypercyclic for\,} \mathsf{T}_a\}$ 
	must have Hausdorff dimension zero.
	\qed
\end{thm}

We can extend Theorem~\ref{zero measure 1 dim} as follows.

\medskip\noindent {\bf Theorem C.}
{\it
For any integer $n\geqslant 1$, $m\geqslant 1$,
 there exists some transcendental growth function $\hat{\phi}$ of the shape~\eqref{growth phi},
	such that for every
	$F\in \mathsf{Hol}(\mathbb{C}^n,\mathbb{C}^m)$ with slow growth
	\begin{equation*}
	\norm{F(z)}\leqslant
		\hat{\phi}(\norm{z})
		\qquad
		\qquad
		{\scriptstyle
			(\forall\,z\,\in\, \mathbb{C}^n),
		}
	\end{equation*}
the set 
\begin{equation}
    \label{set of hypercyclic directions on spheres}
    I_{F}=\{
    a\in  \mathbb{S}^{2n-1}: F \text{\,is hypercyclic for}\, \mathsf{T}_a
    \}
\end{equation}
	must have zero Lebesgue measure.
}

\medskip
For $n=1$, there holds a finer result,
see Observation~\ref{observation for n=1}.

\begin{ques}
Denote $\mathcal{T}=\{\phi:\mathbb{R}_{\geqslant 0} \rightarrow \mathbb{R}_{> 0}\,\,\, \text{is a continuous function satisfying~\eqref{growth phi}}\}$. Given positive integers $n, m\geqslant 1$, 
what is the value of \[
\inf_{\phi\in \mathcal{T}}\,\big(\sup\,\{\dim_{
\mathscr{H}} I_F :  F\in \mathsf{Hol}(\mathbb{C}^n,\mathbb{C}^m)\cap\mathsf{S}_\phi\}\big)
? 
\]
Here $\dim_{\mathscr{H}}$ stands for the Hausdorff dimension. 
\end{ques}

\begin{ques} 
Let ${\phi}$ be given as~\eqref{growth phi}.  Find some sufficient (resp. necessary) condition on an uncountable direction set $I\subset \mathbb{S}^{2n-1}$ so (resp. such)  that 
$
\cap_{a\in I}\mathsf{HC}(\mathsf{T}_a)\cap\mathsf{S}_\phi\neq\emptyset
$
(compared with Theorem~\ref{main in math ann}).
\end{ques}

\begin{defi}[\cite{MR2554568}]\label{Oka mfd}
 A complex manifold $Y$ is an Oka manifold if every holomorphic map 
from a neighbourhood of a compact convex set $K$ in a Euclidean space $\mathbb{C}^n$ (for any $n\in \mathbb{Z}_+$) to
$Y$ is a uniform limit on $K$ of entire maps $\mathbb{C}^n \rightarrow Y$.
\end{defi}

Known examples of Oka manifolds include $\mathbb{C}^n$,  $\mathbb{CP}^n$, complex tori,  complex
homogeneous manifolds~\cite{MR0098197}, elliptic  manifolds~\cite{MR1001851},    
 hyperquadrics~\cite{MR2350038}, smooth toric varieties~\cite{Lrusson2011SmoothTV},  ``ball complements''~\cite{kusakabe-annals}, etc~\cite{Forstneric-Oka-book, MR4547869}.

\medskip

 Zajac~\cite{MR3513554} gave a full characterization of hypercyclic composition operators acting on  holomorphic functions on Stein manifolds. 
Kusakabe~\cite{MR3639046} considered universal
holomorphic maps from Stein manifolds to Oka manifolds rather than to $\mathbb{C}$.
By the   argument of~\cite{MR3679721}, one can show that
the $G_\delta$-density of common hypercyclic elements in the space  $\mathsf{Hol}(
 \mathbb{C}^n, Y)$ of holomorphic maps from $\mathbb{C}^n$  to any connected Oka manifold $Y$ with respect to  translation operators
$\mathsf{T}_a:
f(\bullet)
\mapsto
f(\bullet+a)
$ for all $a\in\mathbb{C}^n\setminus\{\mathbf{0}\}$ (see Theorem~\ref{BG Theorem}).

Next, we establish the abundance of common frequently hypercyclic elements in   $\mathsf{Hol}(
 \mathbb{C}^n, Y)$  with respect to  
$\mathsf{T}_a:
f(\bullet)
\mapsto
f(\bullet+a)
$ for all $a\in\mathbb{C}^n\setminus\{\mathbf{0}\}$. 
A source of inspiration is~\cite{bayart2016common}, in which Bayart
established 
the $G_\delta$-genericity of   hypercyclic elements for 
high-dimensional families of operators
in the framework of $F$-spaces enjoying linear structure.
Note that most Oka manifolds have no linear structure, while
satisfying various equivalent Runge-type approximation properties
(cf.~\cite{MR2199229, Forstneric-Oka-book}).

\smallskip

\begin{defi}[\cite{book-bayart-matheron, book-linear-chaos}]
  For a subset $A$ of natural integers $\mathbb{N}$, the lower density of $A$ is defined as
  \[
  \underline{\mathrm{dens}}(A)\coloneqq
  \liminf_{N\rightarrow +\infty}\,
  \frac{\#\{n\leqslant N: n\in A\}}{N}.
  \]
  We call
  $x\in \mathsf{HC}(\mathsf{T})$ a {\sl frequently hypercyclic  element} for $\mathsf{T}$, if for any nonempty open subset $U$ of the ambient space, the set $\{n\in\mathbb{N}: \mathsf{T}^{n} x\in U\}$ has positive lower density. We denote
  by $\mathsf{FHC}(\mathsf{T})$ the set of frequently hypercyclic elements for $\mathsf{T}$.
\end{defi}

The notion of lower density goes back at least to 
Voronin Universality Theorem about the Riemann zeta-function~\cite{MR0472727}.

\medskip\noindent {\bf Theorem D.}
{\it
Let $Y$ be a connected Oka manifold, and $n\geqslant 1$ an integer. Then in $\mathsf{Hol}(\mathbb{C}^n,Y)$ with the compact-open topology, 
$\cap_{b\in\mathbb{C}^n\setminus\{{\mathbf{0}}\}}\,
    \mathsf{FHC}(\mathsf{T}_b)$ is a dense subset.
}

\smallskip

This   is a surprising result, at least to the authors. The key ingredient in the proof of Theorem~D is to arrange many disjoint closed
hypercubes densely  in some bounded domains in $\mathbb{C}^n$, subject to a  certain equidistribution criterion, such that the union of these closed hypercubes is polynomially convex. This is already an interesting (and difficult) problem by its own.

\medskip
Note that a  certain configuration of hypercubes also appeared in~\cite{MR3679721}, by which Bayart and Gauthier  generalized  Theorem~\ref{all trans}~\cite{amazing-theorem} in higher dimension as 

\begin{thm}[\cite{MR3679721}]\label{BG Theorem}
Let $Y$ be a connected Oka manifold, and $n\geqslant 1$ an integer. Then in  $\mathsf{Hol}(\mathbb{C}^n, Y)$ with the compact-open topology,
$\cap_{a\in\mathbb{C}^n\setminus\{\mathbf{0}\}}\,
    \mathsf{HC}(\mathsf{T}_a)$ contains a $G_\delta$-dense subset. \qed
\end{thm}

The purpose of using  hypercubes rather than  balls is to simplify the argument of polynomial convexity (see Lemma~\ref{poly convex}). In fact, the use of  unions of hypercubes appeared at the  very beginning of several complex variables, see e.g.~\cite{oka1936fonctions}. 
The original statement of~\cite{MR3679721} only concerns $Y=\mathbb{C}$, while the proof actually works  for any connected Oka manifold $Y$, since 
$\mathsf{Hol}(\mathbb{C}^n,Y)$ is separable (see Remark~\ref{separable}) and that Runge's approximation property holds for any holomorphic map from a neighborhood of a 
polynomially convex set $K$  in 
$\mathbb{C}^n$ to $Y$ (see Remark~\ref{why can bopa}), thanks to the fundamental work of Forstneri\v{c}~\cite{MR2199229} asserting that each Oka manifold $Y$ satisfies 
the following {\sl BOPA property} (\cite[p.~258]{Forstneric-Oka-book}).

\medskip\noindent{\it Every continuous map $f_0:X\rightarrow Y$ from a Stein space $X$ that is holomorphic on a neighborhood of a compact $\mathcal{O}(X)$-convex subset $K\subset X$ can be deformed to a holomorphic map $f_1: X\rightarrow Y$ by a homotopy of maps that are holomorphic near $K$ and arbitrarily  close to $f_0$ on $K$.}

\medskip

For Theorem~D, we have to manipulate hypercubes carefully on every  $\mathsf{Q}(0,R+N)\setminus \overline{\mathsf{Q}}(0,R)$ for some $N>0$ independent of $R\gg 1$, while
for  Theorem~\ref{BG Theorem}, one has more flexibility about placing hypercubes outside $\overline{\mathsf{Q}}(0,R)$ for every $R\gg 1$. 
Therefore,  we cannot directly use the construction of~\cite{MR3679721} for Theorem~D, while  we can provide an alternative proof of Theorem~\ref{BG Theorem}
by  our  hypercubes arrangement, see Appendix~\ref{append a}.

\begin{ques}
In spirit of~\cite[Problem 9.1]{Dinh-Sibony-list} and~\cite[Questions 5.1--5.3]{Guo-Xie-1},
we can ask the corresponding questions for frequently hypercyclic   holomorphic maps
from some source space $X$ with mild symmetry to a connected Oka manifold $Y$ (cf.~\cite[Theorem~1.4]{MR3639046}, ~\cite[Remark~3.3]{Xie-IMRN}).  

\smallskip\noindent
(1).  Find the minimal growth rate among entire (resp. meromorphic) functions in $\mathsf{FHC}(\mathsf{T}_a)$ for $X=\mathbb{C}^n, Y=\mathbb{C}^m$ (resp. $Y=\mathbb{CP}^m$), where
$a\in \mathbb{C}^n\setminus \{0\}$. 
Note that the basic case $X=Y=\mathbb{C}$ of this problem was completely solved by~\cite{MR2579678}.

   \smallskip\noindent
(2).  Determine the  minimal growth rate among Nevanlinna characteristic functions $T_f(r)$ (cf.~\cite{MR4265173}) associated with
frequently hypercyclic (with respect to $\mathsf{T}_1$)
 entire curves $f: \mathbb{C}\rightarrow Y$  in a compact connected Oka manifold $Y$.
 
    \smallskip\noindent
(3). Seek the minimal growth rate of common hypercyclic (resp. frequently hypercyclic) elements for all nontrivial translation operators in the above settings (1) and (2).
\end{ques}

\medskip
Here is the structure of this paper.
In Sect.~\ref{Preliminaries}, we provide some  preparations.
In Sect.~\ref{sect. proof of thm a b} and~\ref{proof of thm D}, we  establish Theorems A+B+C and D respectively.

\bigskip\noindent
{\bf Acknowledgments.}
We thank  Zhangchi Chen and Dinh Tuan Huynh for stimulating conversations.
We thank Yi C. Huang for  nice advice.
We are very grateful to Paul Gauthier for careful reading of the manuscript and for providing excellent suggestions (2 pages!).
We thank the referee for nice suggestions
and for pointing out some important references.

\medskip\noindent
{\bf Funding.}
The second named author  is
partially supported by 
National Key R\&D Program of China Grants
No.~2021YFA1003100, 
No.~2023YFA1010500 and  National Natural Science Foundation of China Grant No.~12288201.

\section{\bf Preparations}\label{Preliminaries} 

\subsection{Hypercyclic operator theory}
We refer the interested readers to the nice  expositions~\cite{book-bayart-matheron,	 book-linear-chaos} for learning this subject. 

\begin{defi}
Let $X$ be a metric space.
	A continuous operator $\mathsf{T}: X\rightarrow X$ is called \textit{topologically transitive}, if for any two nonempty open subsets $U$, $V$ of $X$, there exists some $n\geqslant 1$ such that $\mathsf{T}^{n}(U)\cap V\neq \emptyset$.
\end{defi}

The equivalence between hypercyclicity and topological transitivity
for continuous operators $\mathsf{T}$ on a separable Fr\'echet space 
is known as the {\sl Birkhoff transitivity theorem} 
(cf. e.g.~\cite[p.~10]{book-linear-chaos}). 
Moreover, from the proof, one sees that

\begin{pro}\label{countable direction}
	On a separable Fr\'echet space, any countable family of hypercyclic operators $\{\mathsf{T}_n\}_{n\in\mathbb{N}}$ has a dense set of common hypercyclic elements, {\em i.e.}, $\bigcap_{n\in\mathbb{N}}\mathsf{HC}(\mathsf{T}_n)$ is a $G_{\delta}$-dense subset. \qed
\end{pro}


\begin{defi}
	Let $X$, $Y$ be two metric spaces,
	and $S:Y\rightarrow Y$, $T:X\rightarrow X$ be continuous maps. 
	We say that $T$ is  {\sl quasi-conjugate} to $S$
	if there exists a continuous map $\phi:Y\rightarrow X$  with dense image such that the following diagram commutes
	\[
	\begin{tikzcd}
		Y \arrow[r, "S"] \arrow[d, "\phi"] & Y \arrow[d, "\phi"] \\
		X \arrow[r, "T"]                   & X.                  
	\end{tikzcd}
	\]
\end{defi}

\begin{pro}\label{hypercyclic comparison}(cf.~\cite[Propositions 1.13, 1.19]{book-linear-chaos})
	 The following are preserved under quasi-conjugacy:
\begin{itemize}
	\smallskip
	\item[$\diamondsuit\, 1.$]
	topological transitivity;
	
	\smallskip
	\item[$\diamondsuit\, 2.$]
dense orbit.
	\qed
\end{itemize}
\end{pro}

\begin{defi}\label{def mixing}
	A continuous operator $\mathsf{T}: X\rightarrow X$ is called {\sl mixing} if for any nonempty open subsets $U$, $V$ of $X$, there exists some $N\geqslant 1$ such that $\mathsf{T}^{n}(U)\cap V\neq \emptyset$ for all $n\geqslant N$.
\end{defi}

\subsection{Analytic maps between infinite dimensional complex Euclidean spaces}\label{Infinite dimension analytic function}
For our purpose, 
we desire a separable ``$\mathsf{Hol}(\mathbb{C}^{\infty}, \mathbb{C}^{\infty})$'' which contains all $\mathsf{Hol}(\mathbb{C}^{n}, \mathbb{C}^{m})$ as subspace. However, it is known in the literature~\cite{dineen2012complex} that the space of  ``holomorphic functions'' on an infinite dimensional space cannot be separable. Notice that 
in the finite dimensional case, holomorphic maps $\mathsf{Hol}(\mathbb{C}^{n}, \mathbb{C}^{m})$ are nothing but the analytic maps $\mathcal{A}(\mathbb{C}^{n}, \mathbb{C}^{m})$. Hence, instead of  ``$\mathsf{Hol}(\mathbb{C}^{\infty}, \mathbb{C}^{\infty})$'',
we shall consider some  analytic maps ``$\mathcal{A}(\mathbb{C}^{\infty}, \mathbb{C}^{\infty})$'' 
between infinite dimensional complex spaces.

Denote the set of finite multi-indices by
\begin{equation}
    \label{what is beautiful F}
\mathscr{F}=\{(i_j)_{j\in\mathbb{N}}:\quad\text{every\,} i_j\in\mathbb{Z}_{\geqslant 0},\, \text{all but finite }i_j= 0\},
\end{equation}
which is countable. 
Each 
\[
K=(i_1, \dots, i_n, 0, 0, 0,\dots)\in \mathscr{F}
\]
corresponds to a monomial 
\begin{equation}
    \label{define Z^K}
Z^K\coloneqq z_1^{i_1}\cdots z_n^{i_n}.
\end{equation}
On the infinite dimensional $\mathbb{C}$-linear space
\[
\ell^2(\mathbb{C})=\{(a_i)_{i\in\mathbb{N}}:\,
\sum_{i\geqslant 1}\,|a_i|^2<+\infty\}
\]
equipped with the weak topology, 
we say a continuous function $f:\ell^2(\mathbb{C})\rightarrow\mathbb{C}$ is analytic if 
it can be written as a power series $\sum_{I\in \mathscr{F}} a_I\cdot Z^I,\,a_I\in \mathbb{C}$, which is uniformly absolutely convergent on closed balls $\overline{\mathbb{B}}(0,R)$, for all $R>0$.  
The space of analytic functions on $\ell^2(\mathbb{C})$ is denoted by $\mathcal{A}^{\mathsf{w}}(\ell^2(\mathbb{C}),\mathbb{C})$. 
Also, we define
analytic maps from
$\ell^2(\mathbb{C})$ to $\mathbb{C}^m=\oplus_{d=1}^m \mathbb{C}$
 as
\[
\mathcal{A}^{\mathsf{w}}(\ell^2(\mathbb{C}), \mathbb{C}^m)\,
\coloneqq\,
\oplus_{d=1}^m\, \mathcal{A}^{\mathsf{w}}(\ell^2(\mathbb{C}),\mathbb{C})
\]
by functoriality. 
The indicator $\mathsf{w}$ on $\mathcal{A}^{\mathsf{w}}$ emphasizes that $f:\ell^2(\mathbb{C})\rightarrow\mathbb{C}$ is a continuous function with respect to the weak-topology of $\ell^2(\mathbb{C})$.
We endow $\mathcal{A}^{\mathsf{w}}(\ell^2(\mathbb{C}), \mathbb{C}^m)$ with the weak-compact-open topology, {\em i.e.}, open sets are generated by
\[
\mathcal{U}^{\mathsf{w}}(K,V)\coloneqq\{f\in\mathcal{A}^{\mathsf{w}}(\ell^2(\mathbb{C}), \mathbb{C}^m):f(K)\subset V\},
\]
for all weakly compact subsets $K$ in $\ell^2(\mathbb{C})$ and for all open subsets  $V$ in $ \mathbb{C}^m$.

The reason for considering the weak topology on the source space $\ell^2(\mathbb{C})$ is as follows.

\begin{pro}
$\mathcal{A}^{\mathsf{w}}(\ell^2(\mathbb{C}), \mathbb{C}^m)$
  is a separable and metrizable space. 
\end{pro}


\begin{proof}

Fix an increasing positive  sequence $\{R_n\}_{n\geqslant 1}\nearrow +\infty$.
Define
 \[
 d(f,g):=\sum_{n=1}^{+\infty}\frac{1}{2^n}\min\{\rho_n(f-g),1\}
 \qquad{\scriptstyle(\forall\, f,\,g\,\in\,\mathcal{A}^{\mathsf{w}}(\ell^2(\mathbb{C}), \mathbb{C}^m))},
 \]
 where
 \[
 \rho_n(f):=\sup_{z\in \overline{\mathbb{B}}(0,R_n)} \norm{f(z)}_{\mathbb{C}^m}
 \qquad{\scriptstyle(\forall\, f\,\in\,\mathcal{A}^{\mathsf{w}}(\ell^2(\mathbb{C}), \mathbb{C}^m))}.
 \]
 Denote
\begin{equation*}
N_{f,R,\epsilon}\coloneqq\{h\in\mathcal{A}^{\mathsf{w}}(\ell^2(\mathbb{C}), \mathbb{C}^m)):\sup_{z\in\overline{\mathbb{B}}(0,R)} \norm{h(z)-f(z)}_{\mathbb{C}^m}<\epsilon\}
\qquad{\scriptstyle(\forall\, f\,\in\,\mathcal{A}^{\mathsf{w}}(\ell^2(\mathbb{C}), \mathbb{C}^m);\,\forall\,R,\,\epsilon\,>\,0)}.
\end{equation*}
Since  $\overline{\mathbb{B}}(0,R)$ is weakly compact (cf.~\cite[Theorem 3.17]{brezis2011functional}), the sets $\{N_{f,R_n,\frac{1}{m}}\}_{n,m\in\mathbb{Z}_+}$ form a countable open neighborhood basis for each $f$. Hence, $d(\cdot, \cdot)$
 is a metric 
 on $\mathcal{A}^{\mathsf{w}}(\ell^2(\mathbb{C}), \mathbb{C}^m)$ inducing the same weak-compact-open topology.

The separability of
$\mathcal{A}^{\mathsf{w}}(\ell^2(\mathbb{C}), \mathbb{C}^m)$
 relies essentially on the fact that $\mathbb{C}$ contains some countable and dense subfield, say 
$\mathbb{Q}_{c}:=
\mathbb{Q}(\sqrt{-1})
$.

First of all,  $\mathsf{Hol}(\mathbb{C}^n,\mathbb{C}^m)$ is clearly separable, since there is an obvious countable dense subset
\[
\mathscr{P}_{n,m}
:=
\{
(P_1(z), \dots, P_m(z)):P_i(z)\in\mathbb{Q}_c[z_1, \dots, z_n],\,1\leqslant i \leqslant m
\}.
\]

For $F=(\sum_{I\in \mathscr{F}} a^{[i]}_I\cdot Z^I)_{1\leqslant i \leqslant m}\in\mathcal{A}^\mathsf{w}(\ell^2(\mathbb{C}),\mathbb{C}^m)$ and a given neighborhood $N_{F,R,\epsilon}$, since each component of $F$ is uniformly absolutely convergent on closed balls $\overline{\mathbb{B}}(0,R)$, we can find a finite subset $A\subset\mathscr{F}$ such that
\[
\sup_{z\in\overline{\mathbb{B}}(0,R)}|\sum_{I\in \mathscr{F}\setminus A} a^{[k]}_I\cdot Z^I|<\frac{\epsilon}{2\sqrt{m}},
\qquad
{\scriptstyle( 1\,\leqslant\,k\,\leqslant\,m)}.
\]
Set $G:=(\sum_{I\in A} a^{[i]}_I\cdot Z^I)_{1\leqslant i \leqslant m}$. Then it is clear that
\[
\sup_{z\in\overline{\mathbb{B}}(0,R)}\norm{G-F}_{\mathbb{C}^m}
\leqslant
\sup_{z\in\overline{\mathbb{B}}(0,R)}
(\sum_{k=1}^m|\sum_{I\in \mathscr{F}\setminus A} a^{[k]}_I\cdot Z^I|^2)^{\frac{1}{2}}
<\epsilon/2.
\]
Since $A$ is finite, we can find some $N$ such that $G$ only involves the variables $z_1, \dots, z_N$, {\em i.e.}, we can regard $G \in \mathsf{Hol}(\mathbb{C}^N,\mathbb{C}^m)$.
Hence, 
there is some $P\in\mathscr{P}_{N,m}$ such that
\[
\sup_{z\in\overline{\mathbb{B}}(0,R)}\norm{G-P}_{\mathbb{C}^m}
<\epsilon/2.
\]
Therefore
$P \in N_{F, R, \epsilon}$. Thus $\cup_{1\leqslant i <+\infty} \mathscr{P}_{i,m}$ is a countable and dense subset of $\mathcal{A}^\mathsf{w}(\ell^2(\mathbb{C}),\mathbb{C}^m)$.
\end{proof}

In the space
\begin{equation}
    \label{analytic maps between infinitely complex spaces}
\mathcal{A}^{\mathsf{w}}(\ell^2(\mathbb{C}),\ell^2(\mathbb{C}))\coloneqq
\Big{\{}
(f_i)_{i\in\mathbb{N}}\in\prod_{i\in \mathbb{N}}\mathcal{A}^{\mathsf{w}}(\ell^2(\mathbb{C}),\mathbb{C}): \,\forall\, z\in\mathbb{C}\,,\,(f_i(z))_{i\in\mathbb{N}}\in\ell^2(\mathbb{C})
\Big{\}},
\end{equation}
endowed with the weak-compact-open topology, 
we introduce  $\mathcal{A}_{\mathsf{fin}}^{\mathsf{w}}(\ell^2(\mathbb{C}),\ell^2(\mathbb{C}))$ the closure of ${\cup_{1\leqslant m<+\infty}\mathcal{A}^{\mathsf{w}}(\ell^2(\mathbb{C}),\mathbb{C}^m)}$,  
which is metrizable and separable. 
 
Recall a classical result that any compact linear operator can be approximated by a sequence of finite rank operators, and that any limit of a sequence of finite rank operators must be compact (\cite{MR2129625}). Here, elements in $\mathcal{A}^{\mathsf{w}}(\ell^2(\mathbb{C}),\mathbb{C}^m)$ shall be compared with finite rank operators,
and elements in  $\mathcal{A}_{\mathsf{fin}}^{\mathsf{w}}(\ell^2(\mathbb{C}),\ell^2(\mathbb{C}))$ shall be compared with  compact linear operators. Hence we can show that 
$\mathcal{A}_{\mathsf{fin}}^{\mathsf{w}}(\ell^2(\mathbb{C}),\ell^2(\mathbb{C}))$ is a proper subspace of 
$\mathcal{A}^{\mathsf{w}}(\ell^2(\mathbb{C}),\ell^2(\mathbb{C}))$.

Likewise, we can define $\mathcal{A}_{\mathsf{fin}}(\mathbb{C}^n,\ell^2(\mathbb{C}))$ to be the closure of ${\cup_{1\leqslant m<+\infty}\mathcal{A}(\mathbb{C}^n,\mathbb{C}^m)}$ in $\mathcal{A}(\mathbb{C}^n,\ell^2(\mathbb{C}))$,  which is also separable and metrizable.

Let $\Lambda$, $J$ be two
index sets counting the dimension of the  source space
\begin{equation}
    \label{explain Lambda}
\mathbf{F}=
\begin{cases}
  \mathbb{C}^n,\quad & \text{if}\,\,\Lambda=\{1, 
  \dots, n\}  \medskip
  \\
  \ell^2(\mathbb{C}),  &\text{if}\,\,\Lambda=\mathbb{N}
\end{cases},
\end{equation}
and the  target space
\[
\mathbf{G}=
\begin{cases}
  \mathbb{C}^m,\quad & \text{if}\,\,J=\{1, \dots, m\}
  \medskip
  \\
  \ell^2(\mathbb{C}),  &\text{if}\,\,J=\mathbb{N}
\end{cases}.
\]
As the weak topology equals the strong topology on every finite dimensional space, we have 
\[
\mathcal{A}(\mathbb{C}^n,\mathbb{C}^m)=\mathcal{A}^w(\mathbb{C}^n,\mathbb{C}^m),\qquad \mathcal{A}_{\mathsf{fin}}(\mathbb{C}^n,\ell^2(\mathbb{C}))=\mathcal{A}^w_{\mathsf{fin}}(\mathbb{C}^n,\ell^2(\mathbb{C})).
\]
If the target space is finite dimensional, $\mathcal{A}_{\mathsf{fin}}^{\mathsf{w}}(\ell^2(\mathbb{C}),\mathbb{C}^m)$ is equal to $\mathcal{A}^{\mathsf{w}}(\ell^2(\mathbb{C}),\mathbb{C}^m)$. Therefore, we adopt the following notation
\begin{equation}
    \label{all 4 types of F, G}
\mathcal{A}_{\mathsf{fin}}^{\mathsf{w}}(\mathbf{F},\mathbf{G}):=
\begin{cases}
    \mathcal{A}(\mathbb{C}^n,\mathbb{C}^m)=\mathsf{Hol}(\mathbb{C}^n,\mathbb{C}^m),\qquad & \Lambda=\{1,\dots,n\},\,J=\{1,\dots,m\},\\
    \mathcal{A}_{\mathsf{fin}}(\mathbb{C}^n,\ell^2(\mathbb{C})),\qquad&\Lambda=\{1,\dots,n\},\,J=\mathbb{N},\\
    \mathcal{A}^{\mathsf{w}}(\ell^2(\mathbb{C}),\mathbb{C}^m),\,\qquad& \Lambda=\mathbb{N},\,J=\{1,\dots,m\},\\
    \mathcal{A}_{\mathsf{fin}}^{\mathsf{w}}(\ell^2(\mathbb{C}),\ell^2(\mathbb{C})),\,\qquad& \Lambda=\mathbb{N},\,J=\mathbb{N}.
\end{cases}
\end{equation}
In summary, we have the following
\begin{pro}
\label{the topology}
$\mathcal{A}_{\mathsf{fin}}^{\mathsf{w}}(\mathbf{F}, \mathbf{G})$
  is a separable and metrizable space. 
 \qed
\end{pro}

\section{\bf Proofs of Theorems A, B, C} 
\label{sect. proof of thm a b}
\subsection{A key lemma}

The engine is Lemma~\ref{needed version} below, which is an infinite dimensional version of  Theorem~\ref{key thm of DSW97} due to
Desch, Schappacher and Webb~\cite{desch1997hypercyclic}.

\begin{defi}\label{weighted space}
Let $I$ be a countable index set and $v=(v_i)_{i\in I}$ be a sequence of positive numbers indexed by $I$.    The {\sl weighted $\ell^1(v)$ space} over $\mathbb{C}$ (resp. $\mathbb{R}$) is
	\begin{equation}
	\label{ell^1(v)}
	\ell^1(v)=\Big\{(\alpha_i)_{i\in I}:\,
	\alpha_i\in \mathbb{C} \text{ (resp. $\mathbb{R}$)},
	\, \sum_{i\in I} v_i|\alpha_i|<\infty\Big\}
	\end{equation}
with the norm  
\[
\norm{(\alpha_i)_{i\in I}}
=
\sum_{i \in I} 
v_i|\alpha_i|.
\]
\end{defi}

For $I=\mathbb{N}$,
define the backward shift operator 
\[
	B:\ell^1(v)\rightarrow \ell^1(v), \qquad (\alpha_{i})_{i\in \mathbb{N}}\mapsto (\alpha_{i+1})_{i\in \mathbb{N}}.
	\]

\begin{thm}[\cite{desch1997hypercyclic}]\label{key thm of DSW97}
	If $\sup_{i\in\mathbb{N}}\frac{v_i}{v_{i+1}}<+\infty$, then  $
	e^B
	:=
	\sum_{k=0}^{+\infty}\frac{B^k}{k!}
	$ is a mixing operator on $\ell^1(v)$.
\qed
\end{thm}

The original statement in~\cite[Theorem~5.2]{desch1997hypercyclic} only concerns the hypercyclicity of $e^B$, while the proof actually gives the ``mixing'' (see Definition~\ref{def mixing}), cf.~\cite[Theorem 8.1]{book-linear-chaos}. 

\medskip
Let $\Lambda$ be a finite set $\{1, 2, \dots, n\}$ or $\mathbb{N}$ (see~\eqref{explain Lambda}).
Set the multi-index set $\mathscr{I}^{\Lambda}_c$ associated with $\Lambda$ as
\[
\mathscr{I}^{\Lambda}_c=
\begin{cases}
    \mathbb{Z}_{\geqslant 0}^n,\qquad&\text{if}\,\, \Lambda=\{1,\dots,n\},\\
    \mathscr{F}, \qquad&\text{if}\,\, \Lambda=\mathbb{N},
\end{cases}
\]
where
\[
\mathscr{F}=\{(i_j)_{j\in\mathbb{N}}:\quad\text{every\,} i_j\in\mathbb{Z}_{\geqslant 0},\, \text{all but finitely many }i_j= 0\}.
\]
In short, we write multi-indices $(k_i)_{i\in\Lambda}\in \mathscr{I}^{\Lambda}_c$ as $K\in \mathscr{I}^{\Lambda}_c$.

\begin{defi}\label{backward shift}
Let $v=(v_K)_{K\in \mathscr{I}^{\Lambda}_c}$ be a positive weight. For any $s\geqslant 1$, and $$K=(k_1, \dots, k_{s-1}, k_s,  k_{s+1},\dots),$$  define
\[
K_s^{\pm 1}:= (k_1, \dots, k_{s-1}, k_s \pm 1,  k_{s+1},\dots).
\]
We call 
	\[
	B_s:\ell^1(v)\rightarrow \ell^1(v), \qquad (\alpha_K)_{K\in \mathscr{I}^{\Lambda}_c}\mapsto (\alpha_{K_s^{+1}})_{K\in \mathscr{I}^{\Lambda}_c}
	\]
	\textsl{the backward shift operator with respect to the $s$-th index}.
\end{defi}

Introduce a partial order ``$\geqslant$'' on $\mathscr{I}^{\Lambda}_c$ by declaring $K\geqslant K^\prime$ if and only if each index factor of $K$ is no less than the corresponding one of $K^\prime$.
We demand
\begin{equation}\label{continuous condition}
\sup_{{K\geqslant K^\prime}}
\frac{v_{K^\prime}}{v_K}\leqslant 1
\end{equation} 
to ensure that for any $a=(a_i)_{i\in\Lambda}\in \mathbb{C}_c^{\Lambda}$ (recall: $\mathbb{C}_c^{\Lambda}$ is the subset of $\mathbb{C}^\Lambda$ consisting of elements$(x_i)_{i\geqslant 1}$ such that all but finitely many $x_i=0$),
 the finite sum $\sum_{i\in\Lambda} a_i B_i$ and its `exponential'
\[
e^{\sum_{i\in\Lambda} a_i B_i}
=\sum_{k=0}^{+\infty}\frac{{(\sum_{i\in\Lambda} a_i B_i)}^k}{k!}
\]
is
continuous on $\ell^1(v)$. Since
\begin{align*}
\norm{B_s((\alpha_K)_{K\in \mathscr{I}^{\Lambda}_c})}_{\ell^1(v)}
=&
\sum_{K\in\mathscr{I}^{\Lambda}_c} v_{K}|\alpha_{K_s^{+1}}|\\
[\,\text{use \eqref{continuous condition}}]\qquad\leqslant&
\sum_{K\in\mathscr{I}^{\Lambda}_c} v_{K_s^{+1}}|\alpha_{K_s^{+1}}|\\
\leqslant&
\sum_{K\in\mathscr{I}^{\Lambda}_c} v_{K}|\alpha_{K}|=\norm{(\alpha_K)_{K\in \mathscr{I}^{\Lambda}_c}}_{\ell^1(v)},
\end{align*}
we conclude $\norm{B_s}_{\ell^1(v)\rightarrow\ell^1(v)}\leqslant 1$. Hence
\[
\norm{e^{\sum_{i\in\Lambda} a_i B_i}}_{\ell^1(v)\rightarrow\ell^1(v)}
\leqslant\prod_{i\in\Lambda}
\sum_{k=0}^{+\infty}\frac{{\norm{a_i B_i}^k_{\ell^1(v)\rightarrow\ell^1(v)}}}{k!}\leqslant \prod_{i\in\Lambda}
\sum_{k=0}^{+\infty}\frac{{|a_i|}^k}{k!}=e^{\sum_{i\in\Lambda} |a_i|}<+\infty.
\]


\begin{lem}\label{needed version}
	For any positive weight $v=(v_K)_{K\in \mathscr{I}^{\Lambda}_c}$ satisfying \eqref{continuous condition}, for any nonempty subset $J\subset \mathbb{N}$,  and for any $(a_i)_{i\in\Lambda}\in\mathbb{C}_c^\Lambda\setminus\{\mathbf{0}\}$,  the operator $T=
	\bigoplus_J e^{\sum_{i\in\Lambda}a_i B_i}$   is topologically transitive on $\big(\bigoplus_J\ell^1(v)\big)_{\ell^2}$ (see~\eqref{infinite sum l^p}).
\end{lem}

To prove this lemma, we need the following 

\begin{pro}[\cite{book-linear-chaos}--Proposition 2.41]\label{total mixing} For every $n\geqslant 1$,
 let $\mathsf{T}_n$ be a continuous operator on a separable Banach space $X_n$, with $\sup_{n\in\mathbb{N}}\norm{\mathsf{T}_n}<\infty$. Let $1\leqslant p<\infty$.
 Define
\begin{equation}
    \label{infinite sum l^p}
\Big(\bigoplus_{i=1}^{+\infty} X_i\Big)_{\ell^p}:=
\Big{\{}
X=(x_1, x_2, x_3, \dots)
\,\,:\,\,
x_i\in X_i,\,
\norm{X}_{\ell^p}:=\Big(\sum_{i\geqslant 1} \norm{x_i}^p_{X_i}\Big)^{\frac{1}{p}}<\infty\}. 
\end{equation}
    Then $\bigoplus_{n=1}^\infty \mathsf{T}_n$ is mixing on $(\bigoplus_{n=1}^\infty X_n)_{\ell^p}$ if and only if each operator $\mathsf{T}_n$ is mixing.
    \qed
\end{pro}
\begin{proof}[Proof of Lemma~\ref{needed version}]
First, we consider a  basic case $a=(1, 0, 0, \dots)$.
    Rewrite
\begin{equation}
    \label{decompose L1}
\ell^1(v)=\big(\bigoplus_{(k_2,\dots)\in\mathscr{I}^{\Lambda\setminus\{1\}}_c}\ell^1(v_{\cdot,k_2,\dots})\big)_{\ell^1}.
\end{equation}
By Theorem~\ref{key thm of DSW97}, $e^{B_1}$ is mixing on each factor $\ell^1(v_{\cdot,k_2,\dots})$ in~\eqref{decompose L1}.
Hence by Proposition~\ref{total mixing},  $e^{B_1}$ is mixing on the total space $\ell^1(v)$. 
    Using 
    Proposition~\ref{total mixing} again, we see that
    the operator $\bigoplus_J e^{B_1}$ is mixing on $\big(\bigoplus_J\ell^1(v)\big)_{\ell^2}$.
    By symmetry, for other basic cases $a=(0, \dots, 0, 1, 0, 0, \dots)$, the same statement also holds  true.

Next, we deal with  general cases $a=(a_i)_{i\in\Lambda}\in\mathbb{C}_c^\Lambda\backslash\{\mathbf{0}\}$. Without loss of generality, we can assume that $a_1\neq 0$.   
Our strategy is to construct an auxiliary space $\ell^1(w)$ with positive weight $w=(w_K)_{K\in \mathscr{I}^{\Lambda}_c}$ satisfying 	
\[
\sup_{K\in \mathscr{I}^{\Lambda}_c}
\frac{w_K}{w_{K_1^{+1}}}<1,
\]
and a continuous map $F:\ell^1(w)\rightarrow \ell^1(v) $ with dense image such that the following diagram commutes
\[
\begin{tikzcd}
	\ell^1(w) \arrow[rr, "e^{B_1}"] \arrow[d, "F"]           &  & \ell^1(w) \arrow[d, "F"] \\
		\ell^1(v) \arrow[rr, "e^{\sum_{i\in\Lambda}a_i B_i}"] &  & \ell^1(v).               
\end{tikzcd}
\]

Then by taking the direct sum (recall~\eqref{all 4 types of F, G}), we obtain that $$\bigoplus_J\,e^{\sum_{i\in\Lambda}a_i B_i}\,:\,\big(\bigoplus_J\ell^1(v)\big)_{\ell^2}\longrightarrow \big(\bigoplus_J\ell^1(v)\big)_{\ell^2}$$
is quasi-conjugate to $$\bigoplus_J\,
e^{B_1}:\,\big(\bigoplus_J\ell^1(w)\big)_{\ell^2}\longrightarrow \big(\bigoplus_J\ell^1(w)\big)_{\ell^2}.$$
By the preceding argument, $\bigoplus_J e^{B_1}$ is mixing on $\big(\bigoplus_J\ell^1(w)\big)_{\ell^2}$. Hence
by Proposition~\ref{hypercyclic comparison}, we conclude that $$\bigoplus_J\,
e^{\sum_{i\in\Lambda}a_i B_i}:\,\big(\bigoplus_J\ell^1(v)\big)_{\ell^2}\longrightarrow \big(\bigoplus_J\ell^1(v)\big)_{\ell^2}$$ 
is topologically transitive.

Now, we construct such $w=(w_K)_{K\in \mathscr{I}^{\Lambda}_c}$ and $F:\ell^1(w)\rightarrow\ell^1(v)$.

The backward shift operator $B_s$ can be extended to the larger space $\mathbb{C}^{\mathscr{I}^{\Lambda}_c}\supset\ell^1(v)$ naturally 
\begin{equation*}
B_s:\mathbb{C}^{\mathscr{I}^{\Lambda}_c}\rightarrow \mathbb{C}^{\mathscr{I}^{\Lambda}_c}, \qquad (\alpha_{K})_{K\in \mathscr{I}^{\Lambda}_c}\mapsto (\alpha_{K_s^{+1}})_{K\in \mathscr{I}^{\Lambda}_c}. 
\end{equation*}
Define a linear isomorphism
\begin{equation}\label{def phi_z}
	\varphi_z:
	\mathbb{C}^{\mathscr{I}^{\Lambda}_c}
	\rightarrow
	\mathbb{C}[[z_i:i\in\Lambda]],
	\qquad
	(\alpha_K)_{K\in \mathscr{I}^{\Lambda}_c}\mapsto
	\sum_{K\in \mathscr{I}^{\Lambda}_c} \alpha_K\, \frac{Z^K}{K!},    
\end{equation}
where $Z^K$ is defined as~\eqref{define Z^K} and where $K!$ is understood to be the product of the factorials $i_1!
\cdots i_n!$ for $K=(i_1, \dots, i_n, 0, 0,\dots)$, with the convention that $\frac{Z^K}{K!}=0$ if some $i_j<0$. 
Computing
\begin{align*}
\varphi_z\Big(B_s\big((\alpha_K)_{K\in\mathscr{I}^{\Lambda}_c}\big)\Big)
=&
\sum_{K\in \mathscr{I}^{\Lambda}_c} \alpha_{K_s^{+1}}\, \frac{Z^K}{K!}\\
=&
\sum_{K\in \mathscr{I}^{\Lambda}_c} \alpha_{K}\, \frac{Z^{K_s^{-1}}}{K_s^{-1}!}\\
=&
\frac{\partial}{\partial z_s}
\big(
\sum_{K\in \mathscr{I}^{\Lambda}_c} \alpha_{K}\, \frac{Z^{K}}{K!}
\big)=\frac{\partial}{\partial z_s} 
\Big(
\varphi_z\big(
(\alpha_K)_{K\in\mathscr{I}^{\Lambda}_c}
\big)
\Big),
\end{align*}
we obtain

\begin{obs}\label{key obs}
The 
 diagram
\begin{equation}\label{diagram-1}
\begin{tikzcd}
\mathbb{C}^{\mathscr{I}^{\Lambda}_c} \arrow[rr, "{{\sum_{i\in\Lambda} a_i B_i}}"] \arrow[d, "\varphi_z"] &  & \mathbb{C}^{\mathscr{I}^{\Lambda}_c} \arrow[d, "\varphi_z"] \\
{{\mathbb{C}[[z_i:i\in\Lambda]]}} \arrow[rr, "{{\sum_{i\in\Lambda} a_i \frac{\partial}{\partial z_i}}}"] &  & {{\mathbb{C}[[z_i:i\in\Lambda]]}}.
\end{tikzcd}
\end{equation}
is commutative. \qed
\end{obs}


Changing coordinate system linearly by the law
$$
\begin{cases}
	z_1 
	=
	\phi_1(u) = a_1 u_1,
	\medskip
	\\
	z_i \,
	=
	\phi_i(u) 
	=
	a_iu_1 +u_i,\qquad \text{if}\,\,i\neq 1,
\end{cases}
$$
 $\phi=(\phi_1, \phi_2, \phi_3, \dots): \mathbb{C}^{\Lambda}\rightarrow \mathbb{C}^{\Lambda}$ induces an isomorphism  between coordinate rings
\begin{align*}
    \phi^*: \mathbb{C}[[z_i:i\in\Lambda]]&\rightarrow \mathbb{C}[[u_i:i\in\Lambda]]\\
		f(z)&\mapsto f\circ \phi(u).
\end{align*}
By the chain rule, we have 
$$  \frac{\partial}{\partial u_1}=\sum_{i\in\Lambda}\frac{\partial z_i}{\partial u_1}\frac{\partial}{\partial z_i}=\sum_{i\in\Lambda} a_i \frac{\partial}{\partial z_i}.
$$
This combining with \eqref{diagram-1}, we have the following commutative diagram
	\[
	\begin{tikzcd}
		\mathbb{C}^{\mathscr{I}^{\Lambda}_c} \arrow[r, "\varphi_z"] \arrow[d, "\sum_{i\in\Lambda} a_i B_i" ] & {\mathbb{C}[[z_i:i\in\Lambda]]} \arrow[d, "\sum_{i\in\Lambda} a_i \frac{\partial}{\partial z_i}" ] \arrow[r, "\phi^*"] & {\mathbb{C}[[u_i:i\in\Lambda]]} \arrow[d, "\frac{\partial}{\partial u_1}"] & \mathbb{C}^{\mathscr{I}^{\Lambda}_c} \arrow[l, "\varphi_u"'] \arrow[d, "B_1"] \\
		\mathbb{C}^{\mathscr{I}^{\Lambda}_c} \arrow[r, "\varphi_z"]                                            & {\mathbb{C}[[z_i:i\in\Lambda]]} \arrow[r, "\phi^*"]                                                          & {\mathbb{C}[[u_i:i\in\Lambda]]}                         & \mathbb{C}^{\mathscr{I}^{\Lambda}_c}. \arrow[l, "\varphi_u"']                 
	\end{tikzcd}
	\]

Fix the invertible linear map $F=\varphi_z^{-1}\circ \phi^{*-1}\circ \varphi_u :\mathbb{C}^{\mathscr{I}^{\Lambda}_c}\rightarrow \mathbb{C}^{\mathscr{I}^{\Lambda}_c}$. Using the next lemma, we can choose $w=(w_K)_{K\in I}$ 
	such that, by restriction,  $F:\ell^1(w)\rightarrow\ell^1(v)$ is continuous.
	\begin{lem}\label{construct quasiconjugate}
	Let $\mathbb{F}$ be $\mathbb{C}$ or $\mathbb{R}$. 
	Let $G:\mathbb{F}^I\rightarrow \mathbb{F}^I$ be a linear map such that for each  unit coordinate vector $e_i=(\delta_{ij})_{j\in I}$, $G(e_i)=\sum_{j\in I}a^j_i e_j$ is a finite sum, where $\delta_{i j}$ is the Kronecker symbol. Then for any positive weight $v=(v_i)_{i\in I}$, one can find some positive weight $w=(w_j)_{j\in I}$ such that 
	$G:\ell^1(w)\rightarrow\ell^1(v)$ is a continuous map with operator norm $\norm{G}_{\ell^1(w)\rightarrow\ell^1(v)}\leqslant 1$.
	\end{lem}
	\begin{proof}
	For each $i\in I$,  write $G(e_i)=\sum_{j\in J_i}a^j_i\, e_j$, where $J_i$ is a finite index set. Choose 
	\begin{equation}
	    \label{how to choose w}
	w_i>\sum_{j\in J_i}\,|a_i^j|\,v_j.
	\end{equation}
	For any $\alpha=\sum_{i\in I}\alpha_i\, e_i\in \ell^1(w)$, computing
	     $$G(\alpha)= \sum_{i\in I}\,
	     \alpha_i\,
	     G( e_i)
	     =
	     \sum_{i\in I}\,
	     \alpha_i\,
	     \sum_{j\in J_i}\,
	     a^j_i\, e_j,$$
	     we conclude that
	     \[
	     \norm{G(\alpha)}_{\ell^1(v)}\leqslant
	     \sum_{i\in I}|\,
	     \alpha_i|\,
	     \sum_{j\in J_i}\,
	     |a^j_i|\, v_j
	     \leqslant
	     \sum_{i\in I}|\,\alpha_i|\,w_i
	     =
	     \norm{\alpha}_{\ell^1(w)}.
	     \]
\end{proof}

Besides the requirement~\eqref{how to choose w}, we also demand that $(w_K)_{K\in \mathscr{I}^{\Lambda}_c}$
is increasing with respect to the first index (this can be done by selecting $w_K$ inductively in lexicographical order, or Lemma~\ref{Obs 3.1}). Whence 
\[
\sup_{K\in \mathscr{I}^{\Lambda}_c}
\frac{w_K}{w_{K_1^{+1}}}<\infty,
\]
so $e^{B_1}$ is a well-defined continuous linear operator on $\ell^1(w)$. 
It is clear that $F$ has dense image in $\ell^1(v)$.
This concludes the proof.
\end{proof}

\label{proof of Thm A,B}
Recall \eqref{def phi_z}.
For a positive weight $v=(v_{K})_{K\in \mathscr{I}^{\Lambda}_c}$, the image of $\ell^1(v)\subset\mathbb{C}^{\mathbb{N}_c^\Lambda}$ under $\varphi_z$ is 
\[
\varphi_z(\ell^1(v))
=
\Big{\{}
\sum_{K\in \mathscr{I}^{\Lambda}_c}\, \alpha_K\, \frac{Z^K}{K!}\,:\,
\sum_{K\in \mathscr{I}^{\Lambda}_c}\, |\alpha_K|\,v_K <+\infty
\Big{\}}.
\]

\subsection{Controlling growth rates}

 If $(v_K)_{K\in \mathscr{I}^{\Lambda}_c}$ grows rapidly to some order with respect to $K$, then any sequence $(\alpha_{K})_{K\in \mathscr{I}^{\Lambda}_c}$ in $\ell^1(v)$ shall shrink reasonably fast.   
 Hence the corresponding power series $\sum_{K\in\mathscr{I}^{\Lambda}_c} \alpha_K\, \frac{Z^K}{K!}$ could converge,
 and would have slow growth rate. This insight can be achieved by careful manipulation of a standard technique in
  hypercyclic operator theory, cf.~\cite{Chan-Shapiro-91, desch1997hypercyclic, book-bayart-matheron, gomez2010slow,  book-linear-chaos}.

\medskip
For convenience later, we desire the following

\begin{lem}\label{Obs 3.1}
There is an enumeration  $\iota: \mathbb{N} \overset{\sim}{\rightarrow} \mathscr{I}^{\Lambda}_c$ such that, if $K\geqslant K^\prime$, then
 \begin{equation}
     \label{iota has nice order}
     \iota^{-1}(K)\geqslant \iota^{-1}(K^\prime).
 \end{equation}
 
\end{lem}

\begin{proof}
When $\Lambda=\{1, \dots, n\}$, 
such an enumeration can  easily  be obtained, e.g. by using {\sl degree-lexicographic order} (cf. e.g.~\cite[p.~48]{MR1322960}). 
However, when $\Lambda=\mathbb{N}$, this requires much more  effort.
 Fortunately, we discover an enumeration on $\mathscr{F}$ (see~\eqref{what is beautiful F}) encoded in the structure of natural numbers as follows. 
 
 Let $p_i$ be the $i$-th prime number for every $i\geqslant 1$.  Each positive integer $M$ can be factorized as $M=p_1^{k_{M, 1}}\cdots p_n^{k_{M, n}}$ for some $n\geqslant 1$ and $k_{M, 1}, \dots, k_{M, n}\geqslant 0$.
Thus we can define a map $\iota: \mathbb{N}\rightarrow \mathscr{F}$  sending each $M$ to $(k_{M, 1}, \dots, k_{M, n}, 0, 0, \dots)$.
By the existence and uniqueness of prime decomposition,
$\iota$ is uniquely defined and is bijective. 
Whence~\eqref{iota has nice order} holds trivially.
\end{proof}

\begin{lem}
\label{suitable weight}
    Let $\phi: \,\mathbb{R}_{\geqslant 0}\, \rightarrow\, \mathbb{R}_{> 0}$ be a continuous function satisfying~\eqref{growth phi}. Fix a basis
    $\{{P}_K(z)=\frac{Z^K}{K!}\}_{K\in \mathscr{I}^{\Lambda}_c}$ of  $\,\mathbb{C}[[z_i: i\in\Lambda]]$. Then one can choose $v=(v_K)_{K\in \mathscr{I}^{\Lambda}_c}$ such that~\eqref{continuous condition} and the following considerations hold simultaneously.
\begin{itemize}
\smallskip
\item[$\diamondsuit \,1.$]
 The direct sum of  maps (recall~\eqref{all 4 types of F, G}) 
\begin{align*}
    \bigoplus_J\varphi_x:
	\big(\bigoplus_J\ell^1( v)\big)_{\ell^2}
	&\rightarrow
\mathcal{A}_{\mathsf{fin}}^{\mathsf{w}}(\mathbf{F},\mathbf{G})\\
	\Big((\alpha^j_{K_j})_{K_j\in \mathscr{I}^{\Lambda}_c}\Big)_{j\in J}
	&\mapsto
	\Big(\sum_{K_j\in \mathscr{I}^{\Lambda}_c} \alpha^j_{K_j} {P}_{K_j}(x)\Big)_{j\in J}
\end{align*}
is a well-defined continuous linear map. 

\smallskip
\item[$\diamondsuit \, 2.$]
This map
$
    \bigoplus_J\varphi_x:
	\big(\bigoplus_J\ell^1( v)\big)_{\ell^2}
	\rightarrow
	\mathcal{A}_{\mathsf{fin}}^{\mathsf{w}}(\mathbf{F},\mathbf{G})
$
has dense image.

\smallskip
\item[$\diamondsuit \, 3.$]
For any $\alpha\in\big(\bigoplus_J\ell^1(v)\big)_{\ell^2}$, one has
	\[
		\norm{ \bigoplus_J\varphi_x
		(\alpha)(x)}_{G}
		\leqslant 
		\norm{\alpha}_{\big(\bigoplus_J\ell^1(v)\big)_{\ell^2}}
		\cdot
		\phi(\norm{x})
		\qquad
		{\scriptstyle
			(\forall\,x\,\in\, F\,)
		}.
	\]
\end{itemize}
\end{lem}

\begin{proof}
{\bf Step 1.}
Assume for the moment $\mathbf{G}=\mathbb{C}$. With the help of the enumeration $\iota : \mathbb{N}
\overset{\sim}{\rightarrow}
\mathscr{I}^{\Lambda}_c$  defined above, we relabel 
$
\hat P_k(x):= P_{\iota(k)}(x)$
and select
$\hat v_k:= v_{\iota(k)}$
inductively
for $k=1, 2, 3, \dots$ as follows.

Since $\phi$ grows faster than any polynomial, we can find a sequence of positive radii $\{R_k\}_{k\in\mathbb{N}}$ such that
\begin{equation}\label{estiamte dim 1-1}
	\phi(\norm{x})
	\geqslant
	2^k
	\Big|\hat P_k(x)\Big|
	\qquad
	{\scriptstyle
		(\forall\,\norm{x}\,\geqslant\, R_k).
	}
\end{equation}
We now choose an increasing sequence $(\hat v_k)_{k\geqslant 1}$ inductively (hence \eqref{continuous condition} is satisfied, see~Lemma~\ref{Obs 3.1}) obeying the extra law
\begin{equation}\label{estiamte dim 1-2}
	\hat v_k\geqslant 
	\max\Big\{\frac{1}{\phi(0)}
	\sup_{\norm{x}\leqslant R_k} \Big|\hat P_k(x)\Big|,
	\,\frac{1}{2^k}, \hat v_{k-1}
	\Big\}
	\qquad
	{\scriptstyle
		(k\,=\,1,\, 2,\, 3,\, \dots;\,\, \hat{v}_0\,:=\,1)
	}.
\end{equation}

Then for any $\alpha\in\ell^1(v)$, we have
\begin{align*}
	|\varphi_x(\alpha)(x)|
	\leqslant&
	\sum_{m: \norm{x}\leqslant R_m } |\alpha_m|\cdot \Big|\hat P_m(x)\Big|
	+ 
	\sum_{m:R_m < \norm{x}} |\alpha_m|\cdot \Big|\hat P_m(x)\Big|
	\\
	[\text{by}~\eqref{estiamte dim 1-2}, \eqref{estiamte dim 1-1}]\qquad
	\leqslant& 
	\sum_{m: \norm{x}\leqslant R_m} |\alpha_m|\cdot \hat v_m\,{\phi(0)}
	+
	\sum_{m: R_m<\norm{x} } |\alpha_m| \cdot\frac{1}{2^m}\,
	{\phi(\norm{x})}
	\\
	[\text{by}~ \eqref{estiamte dim 1-2}]\qquad	\leqslant& 
	\sum_{m: \norm{x}\leqslant R_m} |\alpha_m|\cdot\hat v_m\,{\phi(0)}
	+
	\sum_{m: R_m<\norm{x}} |\alpha_m|\cdot\hat v_m\cdot
	{\phi(\norm{x})}
	\\
	\leqslant&
	\sum_{m}|\alpha_m|\cdot\hat v_m\cdot \phi(\norm{x})\\
	=&\norm{\alpha}_{\ell^1( v)}\cdot\phi(\norm{x}).
\end{align*}
This ensures that the image of $\varphi$ is contained in $\mathcal{A}^{\mathsf{w}}(\mathbf{F},\mathbb{C})$, since the power series convergences
	\[
	|\varphi_x(\alpha)(x)| \leqslant \norm{\alpha}_{\ell^1(v)}\cdot\phi(\norm{x}) < +\infty
	\qquad
	{\scriptstyle
	(x\,\in\,F)
	}.
	\]
Moreover, we can see that 
$\varphi_x:
	\ell^1(v)
	\rightarrow
	\mathcal{A}^{\mathsf{w}}(\mathbf{F},\mathbb{C})$ is continuous
	by the estimate
\[
    |\varphi_x(\alpha_1)(x)-\varphi_x(\alpha_2)(x)|
    =
    |\varphi_x(\alpha_1-\alpha_2)(x)|
	\leqslant
	\norm{\alpha_1-\alpha_2}_{\ell^1(v)}\cdot\phi(R)
	\quad
	{\scriptstyle
	(\forall\, \alpha_1,\, \alpha_2\,\in\, \ell^1( v),\, \forall\, \norm{x}\,\leqslant\,R)
	}.
\]

Since the set of all polynomials constitutes
a dense subset of $\mathcal{A}^{\mathsf{w}}(\mathbf{F},\mathbb{C})$, we conclude that
 $\varphi_x$ has dense image.
	 
	 \medskip
	 \noindent
{\bf Step 2.} For general $\mathbf{G}$, 
	we choose $v_k$ inductively
	as~ \eqref{estiamte dim 1-1},~\eqref{estiamte dim 1-2}. Likewise, our goal is to show
	 \begin{equation}
	 \label{our goal estimate}
	 \norm{\bigoplus_J\varphi_x(\alpha)(x)}_{\mathbf{G}}\leqslant\norm{\alpha}_{\big(\bigoplus_J\ell^1( v)\big)_{\ell^2}}
		\cdot
		\phi(\norm{x})
		\qquad
		{\scriptstyle
			(\forall\,x\,\in\, F)
		}.
	  \end{equation}
Repeating the same argument as Step 1,
     for each component of 
    $\alpha=
    \Big((\alpha^j_{K_j})_{K_j\in \mathscr{I}^{\Lambda}_c}\Big)_{j\in J}
    $,
    we can obtain an estimate
\begin{equation}\label{estimate in step 1}
    |\sum_{K_j\in \mathscr{I}^{\Lambda}_c}\alpha^j_{K_j}P_{K_j}(x)| \leqslant     
    \norm{(\alpha^j_{K_j})_{K_j\in \mathscr{I}^{\Lambda}_c}}_{\ell^1( v)}\cdot\phi(\norm{x}) 
    \qquad
	{\scriptstyle
	(\forall\,x\,\in\,F,\,\forall\,j\,\in\,J)
	}.   
\end{equation}
Hence the desired estimate~\eqref{our goal estimate} holds
\begin{align*}
    	\norm{\bigoplus_J\varphi_x(\alpha)(x)}_{\mathbf{G}}
    	\,=\,&
    	\Big(\sum_{j\in J} |\sum_{K_j\in \mathscr{I}^{\Lambda}_c}\alpha^j_{K_j}P_{K_j}(x)|^2\Big)^{\frac{1}{2}}\\
    	[\text{by~\eqref{estimate in step 1}}]\qquad
    	\,\leqslant\,&
    	\Big(\sum_{j\in J} \big(\norm{(\alpha^j_{K_j})_{K_j\in \mathscr{I}^{\Lambda}_c}}_{\ell^1(v)}\cdot\phi(\norm{x})\big)^2\Big)^{\frac{1}{2}}\\
    	\,=\,&
    	\Big(\sum_{j\in J} \norm{(\alpha^j_{K_j})_{K_j\in \mathscr{I}^{\Lambda}_c}}^2_{\ell^1(v)}\Big)^{\frac{1}{2}}\cdot\phi(\norm{x})\\
    	\,=\,&\norm{\alpha}_{\big(\bigoplus_J\ell^1(v)\big)_{\ell^2}}
		\cdot
		\phi(\norm{x})
		\qquad
		{\scriptstyle
			(\forall\,x\,\in\, F)
		}.
\end{align*}
Therefore, $
    \bigoplus_J\varphi_x:
	\big(\bigoplus_J\ell^1(v)\big)_{\ell^2}
	\rightarrow
	\mathcal{A}_{\mathsf{fin}}^{\mathsf{w}}(\mathbf{F},\mathbf{G})
$ is
continuous and has dense image.
\end{proof}		

\subsection{A uniform proof for Theorems A, B}
    Take a basis $\{P_K(z)=\frac{Z^K}{K!}\}_{K\in \mathscr{I}^{\Lambda}_c}$  of $\mathbb{C}[z_i:i\in\Lambda]$. 
Choose a positive weight $v=(v_K)_{K\in \mathscr{I}^{\Lambda}_c}$ as in Lemma~\ref{suitable weight}.
Notice that for
$a=(a_1, \dots, a_n, 0, 0, \dots)$,
we have
\begin{equation}\label{trans and exp^diff}
\mathsf{T}_a(\sum_{K} \alpha_K Z^K)
=e^{\sum_{i=1}^{n}a_i \partial_i}(\sum_{K} \alpha_K Z^K) 
\end{equation}
by the Taylor expansion. 
Combining with~\eqref{diagram-1} and Lemma~\ref{suitable weight},	
we obtain the following quasi-conjugate relation for every $a^{[k]}=(a^{[k]}_i)_{i\in\Lambda}\in\mathbb{C}^\Lambda_c\setminus\{\mathbf{0}\}$, $k\geqslant 1$,
\[
\begin{tikzcd}
\big(\bigoplus_J\ell^1(v)\big)_{\ell^2} \arrow[rrr, "\bigoplus_J\,e^{\sum_{i\in\Lambda}a^{[k]}_i B_i}"] \arrow[d, "\bigoplus_J\varphi_z"] &  & & \big(\bigoplus_J\ell^1(v)\big)_{\ell^2} \arrow[d, "\bigoplus_J\varphi_z"] \\
{\mathcal{A}_{\mathsf{fin}}^{\mathsf{w}}(\mathbf{F},\mathbf{G})} \arrow[rrr, "\mathsf{T}_{a^{[k]}}"]                                                                                          &  &  & {\mathcal{A}_{\mathsf{fin}}^{\mathsf{w}}(\mathbf{F},\mathbf{G})}.                                                     
\end{tikzcd}
\]  

Condition \eqref{continuous condition}, which is guaranteed by Lemma~\ref{suitable weight}, ensures that each $e^{\sum_{i\in\Lambda}a^{[k]}_i B_i}$ is continuous. Then Lemma \ref{needed version} applies to show that the operator $\bigoplus_J e^{\sum_{i\in\Lambda}a^{[k]}_i B_i}  : \big(\bigoplus_J\ell^1(v)\big)_{\ell^2}\rightarrow \big(\bigoplus_J\ell^1(v)\big)_{\ell^2}$ is hypercyclic. 
It follows from Proposition~\ref{countable direction} that $\cap_{k\in\mathbb{N}}\mathsf{HC}(\bigoplus_J e^{\sum_{i\in\Lambda}a^{[k]}_i B_i})$ is a $G_\delta$-dense subset  for any countable $\{a^{[k]}\}_{k\in\mathbb{N}}\subset\mathbb{C}_c^\Lambda$.
    Let $\mathbb{B}(0,1)$ be the closed unit ball of $\big(\bigoplus_J\ell^1(v)\big)_{\ell^2}$. 
    For any $\alpha\in\cap_{k\in\mathbb{N}}\mathsf{HC}(\bigoplus_J e^{\sum_{i\in\Lambda}a^{[k]}_i B_i})\cap \mathbb{B}(0,1)$,
    by Lemma~\ref{suitable weight}, we have
    \[
		\norm{\bigoplus_J\varphi_x(\alpha)(z)}_{\mathbf{G}}
		\leqslant 
		\phi(\norm{z})
		\qquad
		{\scriptstyle
			(\forall\,z\,\in\, F)
		}.
	\]
	Since $\bigoplus_J\varphi_x$ is continuous and the orbit of $\alpha$ under each operator $e^{\sum_{i=1}^{n}a^{[k]}_i B_i}$ is dense,  by Proposition~\ref{hypercyclic comparison}, the orbit of $ \bigoplus_J\varphi_x(\alpha)(z)$ under  each operator $\mathsf{T}_{a^{[k]}}$ must be also dense in $\mathcal{A}_{\mathsf{fin}}^{\mathsf{w}}(\mathbf{F},\mathbf{G})$.
 Hence $\bigoplus_J\varphi_x(\alpha)(z)$ is a common hypercyclic element for all $\{\mathsf{T}_{a^{[k]}}\}_{k=1}^{+\infty}$. Summarizing, we conclude
 \[
 (\bigoplus_J\varphi_x)
 \Big(
 \cap_{k\in\mathbb{N}}\mathsf{HC}(\bigoplus_J e^{\sum_{i\in\Lambda}a^{[k]}_i B_i})\cap \mathbb{B}(0,1)
 \Big)
 \subset
 \cap_{k\in\mathbb{N}} \mathsf{HC}(\mathsf{T}_{a^{[k]}})\cap\mathsf{S}_\phi\subset\mathcal{A}_{\mathsf{fin}}^{\mathsf{w}}(\mathbf{F},\mathbf{G}).
 \]
 For the finite dimensional case where $\mathcal{A}_{\mathsf{fin}}^{\mathsf{w}}(\mathbf{F},\mathbf{G})=\mathsf{Hol}(\mathbb{C}^n,\mathbb{C}^m)$, applying Lemma~\ref{amazing-lem}, we  obtain
 \[
 (\bigoplus_J\varphi_x)
 \Big(
 \cap_{k\in\mathbb{N}}\mathsf{HC}(\bigoplus_J\ e^{\sum_{i\in\Lambda}a^{[k]}_i B_i})\cap \mathbb{B}(0,1)
 \Big)
 \subset
 \cap_{r\in\mathbb{R}_+}\cap_{k\in\mathbb{N}} \mathsf{HC}(\mathsf{T}_{r\cdot a^{[k]}})\cap\mathsf{S}_\phi\subset\mathsf{Hol}(\mathbb{C}^n,\mathbb{C}^m).
 \qed
 \]

\subsection{\bf  Proof of Theorem~C}\label{Proof of theorem C}

\begin{obs}
\label{observation for n=1}
For any integer  $m\geqslant 1$,
 there exists some transcendental growth function ${\phi}$ of the shape~\eqref{growth phi},
	such that for every
	$F\in \mathsf{Hol}(\mathbb{C},\mathbb{C}^m)$ with slow growth
	\begin{equation*}
	\norm{F(z)}\leqslant
		{\phi}(\norm{z})
		\qquad
		\qquad
		{\scriptstyle
			(\forall\,z\,\in\, \mathbb{C}),
		}
	\end{equation*}
the set $$\{a\in  \mathbb{S}^{1}: F \text{\,is hypercyclic for}\, \mathsf{T}_a\}$$  
	has zero Hausdorff dimension.
\end{obs}

\begin{proof}
Take one $\phi$ satisfying Theorem~\ref{zero measure 1 dim}. 
Let $F_1$ be the first coordinate factor of $F=(F_1, \dots, F_m)$.
Note that if $a\in  \mathbb{S}^{1}$ is a hypercyclic direction for $F$,
it must be also for $F_1\in \mathcal{H}(\mathbb{C})$.
By Theorem~\ref{zero measure 1 dim}, we conclude the proof.
\end{proof}

\begin{proof}[Proof of Theorem C] 
We only need to verify the case that $n\geqslant 2$.
Take a growth function $\phi$ of the shape~\eqref{growth phi} such that Theorem~\ref{zero measure 1 dim} holds true. 
  Define $\hat{\phi}(r)=\phi(r/\sqrt{n})$.
   We claim that, for any $F\in\mathsf{S}_{\hat\phi}\cap\mathsf{Hol}(\mathbb{C}^n,\mathbb{C}^m)$,  the set \begin{equation}\label{hypercyclic direction set}
       I_F\coloneqq\{e\in  \mathbb{S}^{2n-1}: F \text{\,is hypercyclic for}\, \mathsf{T}_e\}
   \end{equation}
         of hypercyclic directions of $F$ has Lebesgue measure zero. 
         
First, we show that $I_F$ is  Lebesgue measurable. Take a dense countable subset $\{\psi_j\}_{j\geqslant 1}$ of $\mathsf{Hol(\mathbb{C}^n,\mathbb{C}^m)}$ with respect to the compact-open topology. Define
\[
I(j,k)\coloneqq\{e\in  \mathbb{S}^{2n-1}\,:\, 
\exists\, n\in\mathbb{Z}_+\,\text{s.t.}\,\sup_{|z|\leqslant k}|\mathsf{T}_{n\cdot e} F(z)-\psi_j(z)|< 1/k\, \}
	\qquad
		{\scriptstyle
			(\forall\, k,\,j\,\in\,\mathbb{Z}_+)},
\]
which
are clearly open.
Whence $I_F=\cap_{j\geqslant 1}\cap_{k\geqslant 1} I(j,k)$ is measurable.
         
Next, we prove that $I_F$ has zero Lebesgue measure. 
For any $$a=(a_1,\dots,a_n)\in I_F\subset\mathbb{S}^{2n-1}\subset \mathbb{C}^n,$$ 
    find some $|a_j|=\max_{1\leqslant i \leqslant n}|a_i|$. Then define
    \begin{equation}
        \label{define f_a}
    f_a(z)\coloneqq F
    \Big(\frac{a_1}{a_j}z,\dots,\frac{a_{j-1}}{a_j}z,z,\frac{a_{j+1}}{a_j}z,\dots,\frac{a_n}{a_j}z
    \Big),
    \end{equation}
which is clearly an element in $\mathsf{S}_{\phi}\cap \mathsf{Hol}(\mathbb{C},\mathbb{C}^m)$, since $F\in\mathsf{S}_{\hat\phi}$ guarantees that
\[
\norm{f_a(z)}
\leqslant
\hat{\phi}
\Big(
\big(\sum_{i=1}^n\,\big|\frac{a_i}{a_j}\big|^2\,|z|^2
\big)^{1/2}
\Big)
\leqslant
\hat{\phi}(\sqrt{n}\cdot |z|)
=\phi(|z|).
\]

Mimicking~\eqref{hypercyclic direction set}, we denote $I_{f_a}\subset \mathbb{S}^1$ the set of hypercyclic directions for $f_a$. We claim that  $f_a\in\mathsf{HC}(\mathsf{T}_{a_j})\cap \mathsf{Hol}(\mathbb{C},\mathbb{C}^m)$. 
Indeed, direct computation shows
\[
\mathsf{T}_{n\cdot a_j}f_a(z) = 
F
\Big(
\frac{a_1}{a_j}z+n\cdot a_1,\dots,z+n\cdot a_j,\dots,\frac{a_n}{a_j}z+n\cdot a_n
\Big)=
\mathsf{T}_{n\cdot a}\, F
\Big(\frac{a_1}{a_j}z,\dots,z,\dots,\frac{a_n}{a_j}z
\Big).
\]
For any $g\in\mathsf{Hol}(\mathbb{C},\mathbb{C}^m)$, we can regard it as an element $\hat{g}$ in $\mathsf{Hol}(\mathbb{C}^n,\mathbb{C}^m)$ by setting $\hat{g}(z_1,\dots,z_n)=g(z_j)$. For any compact disc $K$ of $\mathbb{C}$ centered at zero, and any  $\epsilon > 0$, since $a\in I_F$, there is some $n\gg 1$ such that
$\mathsf{T}_{n\cdot a} F(\cdot)$
approximates
$\hat{g}(\cdot)$
uniformly on $K\times\cdots\times K$ ($n$ times) having error less than $\epsilon$.
Therefore, for any $z\in K$,
\begin{equation*}
    \label{some computation on K}
\mathsf{T}_{n\cdot a_j}f_a(z) = \mathsf{T}_{n\cdot a} F\Big(\frac{a_1}{a_j}z,\dots,z,\dots,\frac{a_n}{a_j}z\Big)
\approx
\hat{g}\Big(\frac{a_1}{a_j}z,\dots,z,\dots,\frac{a_n}{a_j}z\Big)
=
g(z),
\end{equation*}
{\em i.e.}, $\mathsf{T}_{n\cdot a_j}f_a$ can approximate $g$
on $K$ up to error $\epsilon$.

Summarizing, 
$f_a\in\mathsf{S}_{\phi}\cap \mathsf{Hol}(\mathbb{C},\mathbb{C}^m)\cap \mathsf{HC}(\mathsf{T}_{a_j})$, whence $a_j/|a_j|\in I_{f_a}$ by~Lemma~\ref{amazing-lem}.  
Now we look at the 
the Hopf fibration
$ \mathbb{S}^1\rightarrow\mathbb{S}^{2n-1}\rightarrow\mathbb{CP}^{n-1}$ 
   induced by
    \[
   \mathbb{CP}^{n-1}:=\frac{\mathbb{C}^n\setminus\{0\}}{\mathbb{C}\setminus\{0\}}\cong\frac{\mathbb{S}^{2n-1}}{\mathbb{S}^1}.
   \]
   The $\mathbb{S}^1$ fiber of
   $[a]\in \mathbb{CP}^{n-1}$ is given by the image of
  \begin{align*}
     i_a: \mathbb{S}^1 & \hookrightarrow\mathbb{S}^{2n-1}\\
      e^{i\theta}&\mapsto (a_1e^{i\theta},\dots,a_n e^{i\theta}).
  \end{align*}
Repeating the above argument, 
 if
$(a_1e^{i\theta},\dots,a_n e^{i\theta})\in I_F$,  setting as~\eqref{define f_a}
\begin{align*}
     f_{e^{i\theta} a}(z)
     &\coloneqq 
     F
    \Big(\frac{a_1 e^{i\theta}}{a_j e^{i\theta}}z,\dots,\frac{a_{j-1} e^{i\theta}}{a_j e^{i\theta}}z,z,\frac{a_{j+1} e^{i\theta}}{a_j e^{i\theta}}z,\dots,\frac{a_n e^{i\theta}}{a_j e^{i\theta}}z\Big)
    \\
    &\,=
    F\Big(\frac{a_1}{a_j}z,\dots,\frac{a_{j-1}}{a_j}z,z,\frac{a_{j+1}}{a_j}z,\dots,\frac{a_n}{a_j}z
    \Big)
    =f_a(z),
\end{align*}
then $$\frac{a_j e^{i\theta}}{|a_j e^{i\theta}|}= \frac{a_j }{|a_j|}\cdot e^{i\theta}\in
I_{f_{e^{i\theta} a}}=I_{f_a}.$$
 By Observation~\ref{observation for n=1}, 
$I_{f_a}$ has  Hausdorff dimension zero, hence has Lebesgue measure zero.

Since the same statement holds for the  intersection of $I_F$ with the Hopf fiber $\mathbb{S}^1$ above any $[a]\in \mathbb{CP}^{n-1}$,
by the Fubini theorem, we conclude that $I_F\subset \mathbb{S}^{2n-1}$ must have $0$ Lebesgue measure.
\end{proof}




\section{\bf Proof of Theorem~D}\label{proof of thm D} 
\subsection{Plan}
Regard $\mathbb{C}^n$ as $\mathbb{R}^{2n}$ by reading each $(z_m=x_{2m-1}+\sqrt{-1}\cdot x_{2m})_{1\leqslant m \leqslant n}\in \mathbb{C}^n$ as $(x_{j})_{1\leqslant j \leqslant 2n}\in \mathbb{R}^{2n}$.
Denote the hypercubes
\begin{align*}
\mathsf{Q}(0,R)\,
&
:=\,\{(x_1, \dots,x_{2n}):|x_i|<R,\,i=1,  \dots, 2n\}\,
\subset\,
\mathbb{R}^{2n}\,
=\,
\mathbb{C}^n
\\
\mathsf{Q}(a,R)\,
&
:=\,
a+\mathsf{Q}(0,R)
	\qquad\qquad
		{\scriptstyle
			(\forall\,
			R\,>\,0;\, \,a\,\in\, \mathbb{R}^{2n})
		}.
\end{align*}

We will prove Theorem~D as a corollary of the following Lemma~\ref{amazing-lem} ($\spadesuit  2$) and Theorem~\ref{dense to Oka}.

\begin{lem}\label{amazing-lem}
    Let $Y$ be a connected Oka manifold. For any $\theta\in \partial\mathsf{Q}(0,1)$ and $ b\in\mathbb{R}_+\cdot \theta$, in 
    $\mathsf{Hol}(\mathbb{C}^n,Y)$ with the compact-open topology, one has
    \begin{itemize}
        \item[($\spadesuit  1$).] $\mathsf{HC}(\mathsf{T}_{\theta})=\mathsf{HC}(\mathsf{T}_{b})$\quad\text{(cf.~\cite[Lemma 1.5]{Guo-Xie-1})};
        
        \smallskip
        \item[($\spadesuit  2$).] $\mathsf{FHC}(\mathsf{T}_{\theta})=\mathsf{FHC}(\mathsf{T}_{b})$.
    \end{itemize}
\end{lem}

Note that ($\spadesuit 1$) is  reminiscent of~\cite[Theorem 8]{amazing-theorem}. We will show
($\spadesuit 2$) (see also a similar result  in~\cite{MR2231886}) in Sect.~\ref{proof of lemma 6.2}.

\begin{thm}\label{dense to Oka}
Let $Y$ be a connected Oka manifold. Then
     $\cap_{\theta\in\partial\mathsf{Q}(0,1)}\,\mathsf{FHC}(\mathsf{T}_\theta)$ contains a dense subset of $\mathsf{Hol}(\mathbb{C}^n,Y)$ with respect to the compact-open topology.
\end{thm}

The proof will be reached in Sect.~\ref{proof of 6.3}, by manipulating a certain sophisticated configuration of  closed hypercubes whose union is polynomially convex.

\subsection{Proof of Lemma~\ref{amazing-lem} ($\spadesuit 2$)}
\label{proof of lemma 6.2}
By
 the argument of~\cite[Theorem 8]{amazing-theorem},
 ($\spadesuit 2$) follows 
 with a minor adaptation (cf.~\cite[Lemma 1.5]{Guo-Xie-1}).
 For readers' convenience, we provide the details here.
 (Note that this result fits into a very general property of semigroup of operators,
see~\cite{MR2294487}.)

\medskip
Fix a complete distance $d_Y$ on $Y$. 
Suppose $f\in\mathsf{FHC}(\mathsf{T}_{\theta})$. 
Our goal is to show that, 
for any $g\in\mathsf{Hol}(\mathbb{C}^n,Y)$ and $r>0$, we can find some set $A\subset\mathbb{N}$ with positive lower density such that
\begin{equation}
\label{approx g goal}
\max_{z\in \overline{\mathsf{Q}}(0,r)}d_Y(\mathsf{T}_{m\cdot b}\,f(z), g(z))<
1/r
\qquad
		{\scriptstyle
			(\forall \,m\,\in\, A)
		}.
\end{equation}

Since $g$ is uniformly continuous
on $\overline{\mathsf{Q}}(0, r+1)$, we
can take a positive number $\delta\ll 1$ such that 
\begin{equation}
\label{g uniformly continuous}
d_Y(g(z_1), g(z_2))<1/3r
\qquad
{\scriptstyle
(\forall\,z_1\,,\,z_2\,\in\,\overline{\mathsf{Q}}(0,r+\delta),\,
\text{with}\,\norm{z_1\,-z_2}\,\leqslant \sqrt{n}\delta)}.
\end{equation}
Choose sufficiently many residue classes $\{[c_i]\}_{i=1}^k$ in $\mathbb{R}\cdot b/ \mathbb{Z}\cdot b$ such that
\begin{equation}
\label{condition on ball centers c_k}
\text{	$\{[c_i]\}_{i=1}^k  $ are $\delta$-dense}
\end{equation} 
in the sense that, any vector $c\in \mathbb{R}\cdot b$ can be approximated by a representative $c'_i$
of some  class $[c_i]$ within tiny distance $c\in\mathsf{Q}(c_i',\delta)$.
This can be achieved by requiring  $k\gg 1$  sufficiently large.

Next,  for $i=1, 2, \dots, k$, we subsequently 
select some representatives $\hat{c}_i\in \mathbb{R}\cdot b$ of  $[c_i]$
far away from each other so that the closed hypercubes
$\overline{\mathsf{Q}}(\hat{c}_i, r+\delta)$ are pairwise disjoint.
Therefore $\cup_{i=1}^k \overline{\mathsf{Q}}(\hat{c}_i, r+\delta)$ is polynomially convex (see Lemma~\ref{poly convex}). 

Now we define a holomorphic map $g_2$ 
on a neighborhood of $\cup_{i=1}^k \overline{\mathsf{Q}}(\hat{c}_i, r+\delta)$ by translations of $g$
\begin{equation}
\label{define g_2}
g_2\restriction_{\overline{\mathsf{Q}}(\hat{c}_i, r+\delta)}
(\bullet)
:=
g(\bullet-\hat{c}_i)
\qquad
{\scriptstyle
	(i\,=\,1,\,2,\, \dots,\,k)}.
\end{equation}
Since $Y$ is Oka, by the BOPA property~\cite[p.~258]{Forstneric-Oka-book}, we can ``extend'' $g_2$ to a holomorphic map $\hat{g}$ from a neighborhood of a large hypercube $\overline{\mathsf{Q}}(0, R)$ with
\begin{equation}
\label{define radius zhuhai R}
R>\max_{1\leqslant i\leqslant k}
\{\norm{\hat{c}_i}\}+r+\delta
\end{equation}
to $Y$ within small error (see Remark~\ref{why can bopa})
\begin{equation}
\label{g_2 near g}
d_Y(\hat{g}(z), g_2(z))
<{1}/{3r}
\qquad
{\scriptstyle
	(\forall\, z\,\in\, \cup_{i=1}^k\, \overline{\mathsf{Q}}(\hat{c}_i, r+\delta))}.
\end{equation}

As $f\in\mathsf{FHC}(\mathsf{T}_{\theta})$, there is some   $\{n_j\}_{j=1}^{+\infty}\subset \mathbb{N}$ with positive lower density such that 
\begin{equation}
\label{f near g}
\sup_{z\in \overline{\mathsf{Q}}(0,R)}d_Y(\mathsf{T}_{n_j\cdot \theta}f(z), \hat{g}(z))<{1}/{3r}
\qquad{\scriptstyle
	(\forall\, j\,\geqslant\, 1)}.
\end{equation}
For each $n_j$,
 at least one vector in $\{n_j\cdot \theta+\hat{c}_i\}_{i=1}^k$ is near to $0+\mathbb{Z}\cdot b$ within the uniform distance $\delta$,
because the  segment between $(- n_j-\delta)\cdot \theta$ and $(-n_j+\delta)\cdot \theta$ cannot avoid all vectors in $\{\hat{c}_i+\mathbb{Z}\cdot b\}_{1\leqslant i\leqslant  k}$ by our construction~\eqref{condition on ball centers c_k}. Hence we can rewrite 
\begin{equation}\label{basic arithmetic}
n_j\cdot \theta+\hat{c}_{i_j} =m_j\cdot b+\delta_j'
\end{equation}
for one  $i_j\in \{1, 2, \dots, k\}$, $m_j\in \mathbb{Z}$ and some residue \begin{equation}
    \label{zhuhai 2}
\norm{\delta_j'}<\delta\norm{\theta}\leqslant \sqrt{n}\delta.
\end{equation}

Summarizing~\eqref{g uniformly continuous}, \eqref{define g_2}, \eqref{g_2 near g},
\eqref{f near g},
for any $z$ in $\overline{\mathsf{Q}}(0, r)$, we obtain
\begin{align*}
\mathsf{T}_b^{m_j}f(z)
&=
f(m_j\cdot b+z)
\\
\text{[from \eqref{basic arithmetic}]}\qquad
&=
f(n_j\cdot a+\hat{c}_{i_j}-\delta_j'+z)
\\
\text{[by~\eqref{define radius zhuhai R}, \eqref{f near g}]}
\qquad
&
\approx 
\hat{g}(\hat{c}_{{i_j}}-\delta_j'+z)
\\
\text{[use~\eqref{g_2 near g}]}
\qquad
&
\approx 
g_2(\hat{c}_{i_j}-\delta_j'+z)
\\
\text{[see~\eqref{define g_2}]}
\qquad
&
=
g(-\delta_j'+z)
\\
\text{[check~\eqref{zhuhai 2}, \eqref{g uniformly continuous}]}
\qquad
&
\approx 
g(z).
\end{align*}
Thus the estimate~\eqref{approx g goal} is guaranteed by adding up small differences in the above three ``$\approx$''.

Lastly, 
it is clear that  for sufficiently large $s$ one has $\{m_j\}_{j\geqslant s}\subset \mathbb{N}$.
We claim that this set
has positive lower density.
A delicate place needs some caution, namely,
different $n_j$ and $n_l$ might contribute the same $m_j=m_l$.
However, noting that  $\{\norm{c_i}+\delta\}_{i=1}^k$ is uniformly bounded,
the repetition multiplicity  of each element in the sequence
$\{m_j\}_{j\geqslant s}$
is uniformly bounded from above by some constant $O(1)<\infty$. Thus we conclude the proof by~\eqref{basic arithmetic} and the positive lower density of $\{n_j\}_{j\geqslant 1}$.
\qed

\begin{rmk}\label{why can bopa}
Now we explain how to use the BOPA property in our setting. Let $F: \cup_i \overline{\mathsf{Q}}_i\rightarrow Y$ be a holomorphic map defined near some  disjoint hypercubes $\{\overline{\mathsf{Q}}_i\}_i$.  
By the connectedness of $Y$, we can  extend 
$F$ to a continuous map $\hat{F}$ from $\mathbb{C}^n$ to $Y$. Hence by the BOPA property, we can deform $\hat{F}$ to some holomorphic map
which approximates $F$ very closely on 
$\cup_i \overline{\mathsf{Q}}_i$. 

Here is  some extra detail.
Every hypercube $\overline{\mathsf{Q}}_i$ is contractible. Hence we can first continuously extend  $F$ to the union of some disjoint neighborhoods $V_i\supset \overline{\mathsf{Q}}_i$  such that the image of each $F(\partial V_i)$ is some  point $y_i\in Y$. As $Y$ is connected, we can find a continuous path $\gamma:\mathbb{R}\rightarrow Y$ passing through these points $\{y_i\}_i$, say $\gamma(i)=y_i$ for each enumeration $i$. Now, we define a continuous map $G: \cup_i \partial V_i\rightarrow \mathbb{R}
$
mapping each $\partial V_i$ to the constant $i$.
By Tietze's extension theorem, we can extend $G$
to some continuous map $\hat{G}: \mathbb{C}^n\rightarrow \mathbb{R}$. Thus we can write a continuous extension of $F$:
\[
\hat{F}(z)=\begin{cases}
    F(z), \qquad& \forall\, z\in\cup_i V_i,\\
    \gamma\circ \hat{G}(z), \qquad& \forall\, z\notin\cup_i V_i.
\end{cases}
\]
\end{rmk}

\subsection{Hypercubes arrangement}
\label{Hypercubes arrangement}

\begin{defi}
     For a compact subset $K\subset\mathbb{C}^n$, the polynomially convex hull of $K$ is
\[
\widehat{K}:=\{w\in\mathbb{C}^n:\,|P(w)|\leqslant\sup_{z\in K}|P(z)|,\,\forall P\in\mathbb{C}[z_1,\dots,z_n]\}.
\]
 $K$ is called polynomially convex if  $K=\widehat{K}$. 
\end{defi}


\begin{lem}[\cite{MR3679721}--Lemma 2.3]
\label{poly convex}
   For two polynomially convex compact subsets $X,Y\subset\mathbb{C}^n=\mathbb{R}^{2n}$, if they can be separated by some real affine coordinate  hyperplane $\{x\in\mathbb{R}^{2n}:x_i=a\}$, for some $1\leqslant i \leqslant 2n, a\in \mathbb{R}$, then $X\cup Y$ is polynomially convex. \qed
\end{lem}
The proof of the above  lemma is straightforward. We thank  an anonymous referee for pointing out that such a lemma can be traced back at least to ~\cite{MR0179383} as a corollary of Kallin's lemma. Now we use Lemma~\ref{poly convex} to provide our key construction.

\begin{lem}\label{placing cubes}
Let $0<\delta< 1/2$ and $k\in\mathbb{N}$ be given. Then 
one can find positive integers $N_0=N_0(k,\delta),\,R_0=R_0(k,\delta)$, such that, for any integer $N\geqslant N_0,\, R\geqslant R_0$, the domain $\mathsf{Q}(0,R+N)\setminus\overline{\mathsf{Q}}(0,R)$ contains
a family of pairwise disjoint closed hypercubes $\{\overline{\mathsf{Q}}(a_i,k+\delta)\}_{i\in\Gamma_{R,k,\delta}}$  with the following two properties.
\begin{itemize}
    \item[($\diamondsuit  1$).] For any $\theta \in\partial\mathsf{Q}(0,1)$, for some $m\in\mathbb{N}$, the  hypercube $\mathsf{Q}(m\cdot\theta,k)$ is contained in one of  $\{\overline{\mathsf{Q}}(a_i,k+\delta)\}_{i\in\Gamma_{R,k,\delta}}$.
    
    \smallskip
    \item[($\diamondsuit  2$).] For any $L\in (0, R)$, the union of the disjoint hyercubes 
    \[
    \cup_{i\in\Gamma_{R,k,\delta}} \overline{\mathsf{Q}}(a_i,k+\delta)\cup\overline{\mathsf{Q}}(0,L)
    \]
    is polynomially convex.
\end{itemize}
\end{lem}

\bigskip
For $m=1,\dots,2n$, consider the closed cones
\begin{align*}
\mathsf{F}_{2m-1}
\,
&:=\,{\{(x_1,x_2,\dots,x_{2n}) : x_{m}\geqslant 0;\,|x_j| \leqslant |x_{m}|,\,\forall j\neq m\}},
\\
\mathsf{F}_{2m}\,
&:=\,{\{(x_1,x_2,\dots,x_{2n}) : x_{m}\leqslant 0;\,|x_j| \leqslant |x_{m}|,\,\forall j\neq m\}}.
\end{align*}
We can use these cones to partition $\partial \mathsf{Q}(0,1)$ as the union of $\{\partial \mathsf{Q}(0,1)\cap \mathsf{F}_s\}_{1\leqslant s \leqslant 4n}$. Likewise,
\[
\mathsf{Q}(0,R+N)\setminus\overline{\mathsf{Q}}(0,R)
=
\cup_{1\leqslant s \leqslant 4n}\,\mathsf{F}_s\cap 
\big( \mathsf{Q}(0,R+N)\setminus\overline{\mathsf{Q}}(0,R)\big).
\]
A proof of Lemma~\ref{placing cubes} will be reached after the following more tractable

\begin{lem}\label{placing cubes on 1-face}
Let $0<\delta< 1/2$ and $k\in\mathbb{N}$  be given. Then
one can find positive integers $N=N(k,\delta)$ and $R_0=R_0(k,\delta)$, such that, for any $R\geqslant R_0$, for each $1\leqslant s \leqslant 4n$,
there exists a family of pairwise disjoint closed hypercubes $\{\overline{\mathsf{Q}}(a_i^{[s]},k+\delta)\}_{i\in\Gamma^{[s]}_{R,k,\delta}}$, with centers $a_i^{[s]}$ in $\mathsf{F}_s\cap\big(\mathsf{Q}(0,R+N)\setminus\overline{\mathsf{Q}}(0,R)\big)$, such that the following two properties hold true.
\begin{itemize}
\smallskip
    \item[($\clubsuit  1$).] For any $\theta \in\partial\mathsf{Q}(0,1)\cap\mathsf{F}_s$, 
    for some $m\in\mathbb{N}$,
    the hypercube $\mathsf{Q}(m\cdot\theta,k)$ is  contained in one of $\{\overline{\mathsf{Q}}(a^{[s]}_i,k+\delta)\}_{i\in\Gamma^{[s]}_{R,k,\delta}}$;
    \item[($\clubsuit  2$).] The  union of disjoint hypercubes
    \[
    \cup_{i\in\Gamma^{[s]}_{R,k,\delta}} \overline{\mathsf{Q}}(a^{[s]}_i,k+\delta)
    \]
    is polynomially convex.
\end{itemize}
\end{lem}

\begin{proof}
By symmetry, we can assume $s=1$.
Take the integer $$l:=\lfloor\frac{4k+2}{\delta}+1\rfloor+1$$ 
where $\lfloor\cdot\rfloor$ means rounding down. Select
\[
N:=(2k+1)\cdot l^{2n-1}+1, \quad R_0:= (2k+1)\cdot( l^{2n-1}-2),
\]
\[
r_j:=(2k+1)\cdot j+R
	\qquad
	{\scriptstyle
	(j\,=\,1\,,\,\dots\,,\, l^{2n-1})
	},
	\]
	\[
	A:=\frac{2k+1}{2k+1+R},
\quad
\delta_0:=\frac{\delta}{(2k+1)\cdot l^{2n-1}+R},
\]
\[
\mathscr{A}:=\{
(1,s_2\, A,s_3\, A,\dots,s_{2n}\, A)
:(s_2,\dots,s_{2n})\in\mathbb{Z}^{2n-1}\},
\]
\[
\mathscr{B}:=\{
(0,t_2\, \delta_0,t_3\, \delta_0,\dots,t_{2n}\, \delta_0)
:(t_2,\dots,t_{2n})\in\mathbb{Z}_{\geqslant 0}^{2n-1};\,t_i\leqslant l-1,\,\forall\,i=2,3, \dots,2n\}.
\]
Enumerate $\mathscr{B}$ as $\{b_i\}_{i=1}^{ l^{2n-1}}$.
Consider the following grid points on $\partial\mathsf{Q}(0,1)\cap\mathsf{F}_1$
\begin{equation}
    \label{A_j zhuhai}
\mathscr{A}_{j}:=\big(b_j+\mathscr{A}\big)\cap\mathsf{F}_1
    \qquad
    {\scriptstyle
		(j\,=\,1\,,\,\dots\,,\, l^{2n-1})
	}.
\end{equation}
We claim that the family of hypercubes $$\{\overline{\mathsf{Q}}(a^{[1]}_i,k+\delta)\}_{i\in\Gamma^{[1]}_{R,k,\delta}}:=\cup_{1\leqslant j\leqslant l^{2n-1}}\{
\overline{\mathsf{Q}}\big(r_j\cdot a,k+\delta\big)
\}_{a\in\mathscr{A}_j}$$  satisfies the desired properties~($\clubsuit  1$) and~($\clubsuit  2$). 

First, by the estimate
\[
R<r_j
=
(2k+1)j+R
\leqslant
(2k+1)\, l^{2n-1}+R
<
R+N,
\]
we can check that the centers are in the target region  (see~\eqref{A_j zhuhai})
\[
\{a^{[1]}_i\}_{i\in\Gamma^{[1]}_{R,k,\delta}}=\{r_j\cdot a\}_{a\in\mathscr{A}_j,1\leqslant j\leqslant  l^{2n-1}}\subset\big(\mathsf{Q}(0,R+N)\setminus\overline{\mathsf{Q}}(0,R)\big)\cap\mathsf{F}_1.
\]
Next, to show 
($\clubsuit  1$),
we associate each point
$a\in\cup_{1\leqslant j\leqslant l^{2n-1}}\mathscr{A}_j$ (see~\eqref{A_j zhuhai}) with a smaller hypercube
$\overline{\mathsf{Q}}(a,\delta_0)$. Hence  we can check that
\[
\cup_{1\leqslant j\leqslant l^{2n-1}}\cup_{a\in\mathscr{A}_j}\overline{\mathsf{Q}}(a,\delta_0)\supset\partial\mathsf{Q}(0,1)\cap\mathsf{F}_1
\]
by the width estimate
\[
 l\cdot\delta_0
 >
 \frac{4k+2}{\delta}\cdot\frac{\delta}{(2k+1)\, l^{2n-1}+R}
 =
 \frac{4k+2}{(2k+1)\, l^{2n-1}+R}
 \overset{(*)}{\geqslant}
 \frac{2k+1}{2k+1+R}=A,
\]
where the  inequality $(*)$ is due to $R\geqslant R_0= (2k+1)\cdot( l^{2n-1}-2)$.

Hence, any $\theta \in\partial\mathsf{Q}(0,1)\cap\mathsf{F}_1$ is contained in some $\overline{\mathsf{Q}}(a,\delta_0)$ for $a\in\mathscr{A}_j,\,1\leqslant j\leqslant l^{2n-1}$. Since
\[
r_j\cdot\delta_0=((2k+1)j+R)\cdot\frac{\delta}{(2k+1)\cdot  l^{2n-1}+R}\leqslant \delta ,
\]
we obtain that $r_j\cdot \theta\in\overline{\mathsf{Q}}(r_j\cdot a,\delta)$.
Whence
\[
\mathsf{Q}(r_j\cdot \theta,k)\,
\subset\,
\overline{\mathsf{Q}}(r_j\cdot a,k+\delta).
\]

Lastly, we show that the hypercubes in the family $\{\overline{\mathsf{Q}}(a_i,k+\delta)\}_{i\in\Gamma^{[1]}_{R,k,\delta}}$ are pairwise disjoint, and satisfy ($\clubsuit  2$). Indeed, 
for every $j\in \{1, 2,  \dots, l^{2n-1}\}$, since
\[
r_j\cdot A=((2k+1)\,j+R)\cdot \frac{2k+1}{2k+1+R} \geqslant 2k+1 >2 (k+\delta),
\]
the hypercubes $\{\overline{\mathsf{Q}}\big(r_j\cdot a,k+\delta\big)
\}_{a\in\mathscr{A}_j}$ are pairwise disjoint.  Noting that centers $a$ are placed on the grid $\mathscr{A}_{j}$, by using Lemma~\ref{poly convex} repeatedly, we see that
$
\cup_{a\in\mathscr{A}_j}\,\overline{\mathsf{Q}}\big(r_j\cdot a,k+\delta\big)
$
is polynomially convex.
Moreover, 
for any two distinct $j, j^\prime$ in $\{1, 2,  \dots, l^{2n-1}\}$, since
$$
|r_j-r_{j^\prime}|=
(2k+1)\cdot |j-j'|\geqslant 2k+1 >2(k+\delta),
$$
it is clear that
$\cup_{a\in\mathscr{A}_j}\,\overline{\mathsf{Q}}\big(r_j\cdot a,k+\delta\big)$
and
$\cup_{a\in\mathscr{A}_{j'}}\,\overline{\mathsf{Q}}\big(r_{j'}\cdot a,k+\delta\big)$ are
disjoint. Hence by using Lemma~\ref{poly convex} repeatedly, we see that
\[
\cup_{1\leqslant j\leqslant l^{2n-1}}\cup_{a\in\mathscr{A}_j}\overline{\mathsf{Q}}\big(r_j\cdot a,k+\delta\big)
\]
is polynomially convex. 
\end{proof}

\begin{proof}[Proof of Lemma~\ref{placing cubes}]
Borrow  two positive integers $N^\prime=N^\prime(k,\delta)$ and $R_0=R_0(k,\delta)$ from Lemma~\ref{placing cubes on 1-face}.
Select $N_0=(4n+1)\cdot (N^\prime+2k+2)$. 
Now we claim that, for any $R\geqslant R_0$, $N\geqslant N_0$, there are pairwise disjoint closed hypercubes $\{\overline{\mathsf{Q}}(a_i,k+\delta)\}_{i\in\Gamma_{R,k,\delta}}$ contained in the domain $\mathsf{Q}(0,R+N)\setminus\overline{\mathsf{Q}}(0,R)$, which satisfy ($\diamondsuit  1$), ($\diamondsuit  2$).

Take
\[
R_m=R+(N^\prime+2k+2)\cdot m
\qquad
    {\scriptstyle
		(m\,=\,1\,,\,2,\,\dots,\,4n+1)}.
\]
By Lemma~\ref{placing cubes on 1-face}, for each $1\leqslant m \leqslant 4n$, since
\[
(R_{m+1}-k-1)
-
(R_{m}+k+1)
=
N',
\]
there is a family of pairwise disjoint closed hypercubes $\{\overline{\mathsf{Q}}(a^{[m]}_i,k+\delta)\}_{i\in\Gamma^{[m]}_{R_m+k+1,k,\delta}}$, with centers in $\mathsf{F}_m\cap\big(\mathsf{Q}(0, R_{m+1}-k-1)\setminus\overline{\mathsf{Q}}(0, R_{m}+k+1)\big)$, such that
\begin{itemize}
    \item[($\bullet\, 1$).] for any $\theta\in\partial \mathsf{Q}(0,1)\cap \mathsf{F}_m$, there is some $l\in\mathbb{N}$ and some $i\in\Gamma^{[m]}_{R_m+k+1, k, \delta}$ with
    \[
    \mathsf{Q}(l\cdot \theta,k)\subset \overline{\mathsf{Q}}(a^{[m]}_i,k+\delta);
    \]
    \smallskip
    \item[($\bullet\,  2$).] the union of
    $\{\overline{\mathsf{Q}}(a^{[m]}_i,k+\delta)\}_{i\in\Gamma^{[m]}_{R_m+k+1,k,\delta}}$
    is polynomially convex.
\end{itemize}
\smallskip
Our desired family of hypercubes is nothing but
\[
\{\overline{\mathsf{Q}}(a_i,k+\delta)\}_{i\in\Gamma_{R,k,\delta}}:=\cup_{1\leqslant m \leqslant 4n}\,
\{\overline{\mathsf{Q}}(a^{[m]}_i,k+\delta)\}_{i\in\Gamma^{[m]}_{R_m+k+1,k,\delta}}.
\]

Obviously, 
these hypercubes are pairwise disjoint, since for each $1\leqslant m\leqslant 4n$, the family of hypercubes
\[
\{\overline{\mathsf{Q}}(a^{[m]}_i,k+\delta)\}_{i\in\Gamma^{[m]}_{R_m+k+1,k,\delta}}\subset\mathsf{Q}(0,R_{m+1})\setminus\overline{\mathsf{Q}}(0,R_m)
\] 
are pairwise disjoint, and the domains $\{\mathsf{Q}(0,R_{m+1})\setminus\overline{\mathsf{Q}}(0,R_m)\}_{1\leqslant m\leqslant 4n}$ are also pairwise disjoint. We can check that all the hypercubes $\{\overline{\mathsf{Q}}(a_i,k+\delta)\}_{i\in\Gamma_{R,k,\delta}}$ are lying in the domain $$\mathsf{Q}(0,R_{4n+1})\setminus\overline{\mathsf{Q}}(0,R_1)\subset\mathsf{Q}(0,R+N)\setminus\overline{\mathsf{Q}}(0,R).$$

It is clear that  ($\diamondsuit  1$) can be deduced from ($\bullet\,1$), since
\[
\partial \mathsf{Q}(0,1)
=
\cup_{1\leqslant i \leqslant 4n}\,\partial \mathsf{Q}(0,1)\cap \mathsf{F}_i.
\]
Now we check ($\diamondsuit 2$).

For any $L<R$, the coordinate hyperplane 
$\{x\in\mathbb{R}^{2n}:x_1=L+1\}$ separates $\overline{\mathsf{Q}}(0,L)$ and $\cup_{i\in\Gamma^{[1]}_{R_1+k+1,k,\delta}}\overline{\mathsf{Q}}(a^{[1]}_i,k+\delta)$.
By ($\bullet\,2$) and Lemma~\ref{poly convex}, we see that the  union
\[
\cup_{i\in\Gamma^{[1]}_{R_1+k+1,k,\delta}}\overline{\mathsf{Q}}(a^{[1]}_i,k+\delta)\cup\overline{\mathsf{Q}}(0,L)
\]
is 
polynomially convex. Subsequently, for all $j=1, 2, \dots, 2n$, since
\[
\cup_{1\leqslant s\leqslant 2j-2}\cup_{i\in\Gamma^{[s]}_{R_s+k+1,k,\delta}}\overline{\mathsf{Q}}(a^{[s]}_i,k+\delta)\cup\overline{\mathsf{Q}}(0,L)
\quad\text{and}\quad
\cup_{i\in\Gamma^{[2j-1]}_{R_{2j-1}+k+1,k,\delta}}\overline{\mathsf{Q}}(a^{[2j-1]}_i,k+\delta)
\]
can be separated by $\{x\in\mathbb{R}^{2n}: x_j=R_{2j-1}\}$, while
\[
\cup_{1\leqslant s\leqslant 2j-1}\cup_{i\in\Gamma^{[s]}_{R_s+k+1,k,\delta}}\overline{\mathsf{Q}}(a^{[s]}_i,k+\delta)\cup\overline{\mathsf{Q}}(0, L)
\quad\text{and}\quad\cup_{i\in\Gamma^{[2j]}_{R_{2j}+k+1,k,\delta}}\overline{\mathsf{Q}}(a^{[2j]}_i,k+\delta)
\]
can be separated by $\{x\in\mathbb{R}^{2n}:x_{j}=-R_{2j}\}$,   using ($\bullet\, 2$) and Lemma~\ref{poly convex} repeatedly, eventually, we get the polynomial convexity of
$
\cup_{1\leqslant s\leqslant 4n}\cup_{i\in\Gamma^{[s]}_{R_s+k+1,k,\delta}}\overline{\mathsf{Q}}(a^{[s]}_i,k+\delta)\cup\overline{\mathsf{Q}}(0,L)$.
\end{proof}

\subsection{Proof of Theorem~\ref{dense to Oka}}
\label{proof of 6.3}

First of all, we need the following
\begin{lem}[\cite{book-linear-chaos}--Lemma 9.5]\label{lower density lem} 
There are pairwise disjoint subsets $A(l,\nu)$ of $\mathbb{N}$ for all $l, \nu\in\mathbb{N}$, having positive lower density, such that, for any $n$ in $ A(l,\nu)$ and another different $m\neq n$ in some $A(k,\mu)$, one has $n\geqslant \nu$ and
    $
    |n-m|\geqslant \nu+\mu$.
\qed
\end{lem}

Our goal is to show that, for any given holomorphic map $g: \mathbb{C}^n\rightarrow Y$ and $R\in\mathbb{Z}_+$, there exists a holomorphic map $f\in \cap_{\theta\in\partial{\mathsf{Q}}(0,1)}\mathsf{FHC}(\mathsf{T}_{\theta})$ such that 
\begin{equation}\label{appro-g}
  \sup_{z\in\overline{\mathsf{Q}}(0, R)} d_Y(f(z),g(z))<1/R.
\end{equation}

 \begin{rmk}\label{separable}
 It is well-known that, for any two smooth connected  manifolds $X'$ and $Y'$, the smooth mapping space $\mathcal{C}^\infty(X',Y')$ equipped with the compact-open topology is separable.
 In particular, $\mathcal{C}^\infty(\mathbb{C}^n,Y)$ is separable, so is the closed subset $\mathsf{Hol}(\mathbb{C}^n,Y)\subset \mathcal{C}^\infty(\mathbb{C}^n,Y)$. 
 \end{rmk}
 
Take a countable dense subset $\{\psi_i\}_{i\geqslant 1}$ of $\mathsf{Hol}(\mathbb{C}^n,Y)$. 
For each $j,l\in\mathbb{N}$, by the uniform continuity of $\psi_j$ on the compact closed  hypercube $\overline{\mathsf{Q}}(0,l+1)$, we can find a sufficiently small  positive number $\delta^{[l,j]}\ll 1$ such that
\begin{equation}\label{uniform conti on psi_j}
   d_Y(\psi_j(z_1),\psi_j(z_2))<1/2l
   \qquad
   {\scriptstyle(\forall\,z_1\,,\,z_2\,\in\,\overline{\mathsf{Q}}(0,l+1),\,\text{with}\,z_2\,\in\,\overline{\mathsf{Q}}(z_1\,,\,\delta^{[l,j]})\,)}. 
\end{equation}

Applying Lemma~\ref{placing cubes} upon
 $\delta^{[l,j]}\ll 1$ and $l\in\mathbb{N}$, we obtain  positive integers 
\[
N_0^{[l,j]}=N_0^{[l,j]}(l,\delta^{[l,j]}),\quad R_0^{[l,j]}=R_0^{[l,j]}(l,\delta^{[l,j]})
 \qquad
  {\scriptstyle
	(\,\forall\, l,\,j\,\in\,\mathbb{N})}.
\]
Since there is an enumeration on the index set $\mathbb{N}
\times \mathbb{N}$ (see Lemma~\ref{Obs 3.1}), 
 we can subsequently choose positive integers $N^{[l,j]}\geqslant N_0^{[l,j]}$ for each pair $(l,j)$ along the enumeration, such that $(l,2 N^{[l,j]})\neq (l^\prime,2 N^{[l^\prime,j^\prime]})$ if $(l,j)\neq(l^\prime,j^\prime)$. 
Hence Lemma~\ref{lower density lem} yields a family  $\{A(j,2 N^{[l,j]})\}_{ j,l\in\mathbb{N}}$ of pairwise disjoint subsets of $\mathbb{N}$ having positive lower density.
Set
\[
B(j,2 N^{[l,j]})\,
:=\,
\{p\in A(j,2 N^{[l,j]}) :\,p \geqslant R_0^{[l,j]}+N^{[l,j]}+R\},
\]
which also has positive lower density.

Write all elements of $\bigcup_{l,j\in\mathbb{N}}B(j,2 N^{[l,j]})$ in the increasing order as 
\begin{equation}
    \label{n_k important}
\{n_k\}_{k=1}^{+\infty}.
\end{equation}
Then, for each $n_k\in B(j,2 N^{[l,j]})$,  we define $j_k=j,\,l_k=l$.
Assign each $n_k$ a domain
\[
 \mathsf{A}_k:=\mathsf{Q}(0,n_k+ N^{[l_k,j_k]})\setminus\overline{\mathsf{Q}}(0,n_k- N^{[l_k,j_k]})
 \qquad
 {\scriptstyle
	(\forall \,k\,\geqslant\, 1)}.
\]
By  Lemma~\ref{lower density lem}, we have
\begin{equation}\label{separate domain}
n_k+N^{[l_k,j_k]}\,
\leqslant\,
n_{k+1}
-
\big( 2\,N^{[l_{k+1},j_{k+1}]}+2\,N^{[l_k,j_k]}
\big)\,
+\,N^{[l_k,j_k]}
\,<\,
n_{k+1}-N^{[l_{k+1},j_{k+1}]}
\quad
{\scriptstyle
	(\forall\,k\,\geqslant\, 1)}.
\end{equation}
Hence, $\mathsf{A}_k$ are pairwise disjoint for all $k\geqslant 1$. 

Applying Lemma \ref{placing cubes}, in each $\mathsf{A}_k$, we can insert pairwise disjoint closed hypercubes \[
\mathsf{C}_k
=
\{\overline{\mathsf{Q}}(a_i^{[k]},l_k+\delta^{[l_k,j_k]})\}_{i\in\Lambda_{k}}:=\{\overline{\mathsf{Q}}(a_i,l_k+\delta^{[l_k,j_k]})\}_{i\in\Gamma_{n_k- N^{[l_k,j_k]},l_k,\delta^{[l_k,j_k]}}},\]
such that following two properties hold.
\begin{itemize}
    \item[($\blacklozenge  1$).] For any $\theta \in\partial\mathsf{Q}(0,1)$, for some $m\in\mathbb{N}$, the hypercube $\mathsf{Q}(m\cdot\theta,l)$ is contained in 
    one of $\{\overline{\mathsf{Q}}(a_i^{[k]},l_k+\delta^{[l_k,j_k]})\}_{i\in\Lambda_{k}}$.
    
    \smallskip
    \item[($\blacklozenge  2$).] Set
$
\mathsf{B}_k:=\cup_{i\in\Lambda_{k}}\overline{\mathsf{Q}}(a_i^{[k]},l_k+\delta^{[l_k,j_k]})$. Then
for $k=1$,
$
\mathsf{B}_1\cup\overline{\mathsf{Q}}(0,R)
$ is polynomially convex;
for every $k\geqslant 2$, by~\eqref{separate domain},  the  union 
\[
\mathsf{B}_{k}\cup\overline{\mathsf{Q}}(0,n_{k-1}+N^{[l_{k-1},j_{k-1}]})
\]
is also polynomially convex. 

\end{itemize}

We can define a holomorphic map $g_k$ on a neighborhood of each $\mathsf{B}_k$ by translation
\begin{equation}\label{def of g_k}
    g_k(z)\,
    :=\,
    \psi_{j_k}(z-a_i^{[k]})\qquad
    \qquad{\scriptstyle(\forall\, k\,\geqslant \,1;\, \forall\, i\,\in\,\Lambda_k,\,\forall\,z\,\in\,\overline{\mathsf{Q}}(a_i^{[k]},\, k+\delta^{[l_k,j_k]})}.
\end{equation}
Since $Y$ is Oka, we can inductively use the BOPA property (see Remark~\ref{why can bopa}) to construct a sequence of entire maps $\{f_n\}_{n\geqslant 0}$ as follows.
Fix a sequence of positive numbers $\{\epsilon_k\}_{k\geqslant 1}$ shrinking rapidly, to be specified later.
Starting with $f_0:=g$, by ($\blacklozenge  2$) and the BOPA property, we can find a next holomorphic map $f_{1}: \mathbb{C}^n\rightarrow Y$ such that 
\[
\sup_{z\in\overline{\mathsf{Q}}(0,R)} d_Y(f_{1}(z),f_0(z))<\epsilon_{1},\quad\text{and}\quad
\sup_{z\in\mathsf{B}_{1}} d_Y(f_{1}(z),g_{1}(z))<\epsilon_{1}.
\]
When $f_0,f_1,\dots,f_k$ have been selected for $k\geqslant 1$, again we use ($\blacklozenge  2$) and the BOPA property to obtain a subsequent $f_{k+1}$ with
\[
\sup_{z\in\overline{\mathsf{Q}}(0,n_k+N^{[l_k,j_k]})} d_Y(f_{k+1}(z),f_k(z))<\epsilon_{k+1},\quad\text{and}\quad
\sup_{z\in\mathsf{B}_{k+1}} d_Y(f_{k+1}(z),g_{k+1}(z))<\epsilon_{k+1}.
\]
It is clear that
the limit $f: \mathbb{C}^n\rightarrow Y$ of
$\{f_k\}_{k\geqslant 1}$ is a well-defined holomorphic map
satisfying
\begin{equation}\label{total deviation}
   \sup_{z\in\overline{\mathsf{Q}}(0,R)} d_Y(f(z),g(z))<\sum_{i=1}^{+\infty}\epsilon_i,\quad\text{and}\quad
\sup_{z\in\mathsf{B}_{k}} d_Y(f(z),g_{k}(z))<\sum_{i=k}^{+\infty}\epsilon_i. 
\end{equation}
Now we can subsequently choose $\epsilon_k$ for $k=1, 2, 3, \dots$ shrinking rapidly, such that 
\begin{equation}\label{control deviation}
    \sum_{i=1}^{+\infty}\epsilon_i\,<
    \,
    1/R,\quad\text{and}\quad\sum_{i=k}^{+\infty}\epsilon_i<
    1/{2 l_k}.
\end{equation}

Therefore, the obtained $f$  satisfies \eqref{appro-g}. Moreover, we claim that it approximates each $\psi_j$ frequently under the iterations of the translation operator $\mathsf{T}_{\theta}$ for any $\theta\in\partial\mathsf{Q}(0,1)$. 

Indeed, for each $\psi_j$, 
and for any error bound $\epsilon>0$,
by taking a large $l\in\mathbb{N}$  with ${1}/{l}<\epsilon$, we obtain
 a set $B(j,2N^{[l,j]})$ with positive lower density. Recalling~\eqref{n_k important}, for every  $n_t\in B(j,2N^{[l,j]})$ with $t\geqslant 1$,  by ($\blacklozenge  1$), we can find an integer $m_t$ and one hypercube $\overline{\mathsf{Q}}(a,l+\delta^{[l,j]})$ in the family $\mathsf{C}_t$ such that
\begin{equation}\label{small deviation}
    \mathsf{Q}(m_t\cdot\theta,l)\subset
    \overline{\mathsf{Q}}(a,l+\delta^{[l,j]})
\subset
\mathsf{Q}(0,n_t+N^{[l,j]})\setminus\overline{\mathsf{Q}}(0,n_t-N^{[l,j]}).
\end{equation}
By our previous construction,
it is clear that $(l_t,j_t)=(l,j)$.
Therefore, we can estimate
\begin{align*}
    \sup_{z\in\overline{\mathsf{Q}}(0,l)}d_Y(\mathsf{T}_{m_t\cdot\theta}f(z),\psi_j(z))
    =&
    \sup_{z\in\overline{\mathsf{Q}}(m_t\cdot\theta,l)}d_Y(f(z),\psi_j(z-m_t\cdot\theta))
    \\
    [\text{see~\eqref{small deviation}}]\qquad
    \leqslant&
    \sup_{z\in\overline{\mathsf{Q}}(a,l+\delta^{[l,j]})}d_Y(f(z),\psi_j(z-a))\\
    &
    \qquad
    +
    \sup_{z\in\overline{\mathsf{Q}}(a,l+\delta^{[l,j]})}d_Y(\psi_j(z-m_t\cdot\theta),\psi_j(z-a))
    \\
    [\text{use~\eqref{uniform conti on psi_j}}]\qquad
    \leqslant&
    \sup_{z\in\overline{\mathsf{Q}}(a,l+\delta^{[l,j]})}d_Y(f(z),\psi_j(z-a))+{1}/{2l}
    \\
    [\text{check\,\eqref{def of g_k},\,\eqref{total deviation}}]\qquad
    \leqslant&
    \,\sum_{i=t}^{+\infty}\epsilon_i
    + 
    {1}/{2l}
    \\
    [\text{by\,\eqref{control deviation}}]\qquad
    <& \, {1}/{2l}+{1}/{2l}\,=\, 1/l \,<\, \epsilon.
\end{align*}

Lastly, \eqref{small deviation} implies that $|m_t-n_t|\leqslant N^{[l,j]}$. 
By Lemma~\ref{lower density lem}, any two distinct integers $n_t,\,n_{t^\prime}$ in $B(j,2N^{[l,j]})$ satisfy $|n_{t^\prime}-n_t|\geqslant 4N^{[l,j]}$. 
Therefore
\[
|m_{t'}-m_t|
=
 |(n_{t'}-n_t)-(n_{t'}-m_{t'})
+(n_t-m_t)|
\geqslant
4\,N^{[l,j]}-N^{[l,j]}-N^{[l,j]}
>0.
\]
Hence  the corresponding sequence of integers $\{m_t\}_{n_t\in B(j,2N^{[l,j]})}$
 has positive lower density, inherited from that of
 $B(j,2N^{[l,j]})$.
\qed

\appendix
\section{\bf An alternative proof of Theorem~\ref{BG Theorem}   }\label{append a}
The proof is inspired by that of Theorem~\ref{all trans}~\cite{amazing-theorem} in which $n=1$, $Y=\mathbb{C}$. 

\smallskip
By Lemma~\ref{amazing-lem}~($\spadesuit  1$), we only need to prove the following 

\begin{thm}
Let $Y$ be a connected Oka manifold, and $n\geqslant 1$ an integer. Then in  $\mathsf{Hol}(\mathbb{C}^n, Y)$ with the compact-open topology,
$\cap_{a\in\partial\mathsf{Q}(0,1)}\,
    \mathsf{HC}(\mathsf{T}_a)$ contains a $G_\delta$-dense subset. 
\end{thm}

Take a dense subset $\{\psi_j\}_{j\geqslant 1}$ of $\mathsf{Hol}(\mathbb{C}^n,Y)$ (See Remark~\ref{separable}). 
For any $(a,b,k,j)\in\mathbb{N}^4$, set
\begin{align*}
    E(a,b,k,j)=\{f\in\mathsf{Hol}(\mathbb{C}^n,Y) :\,\forall\, \theta\in\partial\mathsf{Q}(0,1),\exists\, \text{integer}\, s=s(\theta) \leqslant a\,\\
    \text{such that}\,
\sup_{z\in\overline{\mathsf{Q}}(0,k)}d_Y(\mathsf{T}_{s\cdot \theta}f(z),\psi_j(z))<{1}/{b}\}.
\end{align*}
By Definition~\ref{hypercyclic},  $\cap_{\theta\in\partial\mathsf{Q}(0,1)}\mathsf{HC}(\mathsf{T}_\theta)$ contains  $\cap_{b,k,j\geqslant 1}\big(\cup_{a\geqslant 1}E(a,b,k,j)\big)$. Hence, we only need to show that
$\cup_{a\geqslant 1}E(a,b,k,j)$ is a dense open set.

First,
every $E(a,b,k,j)$ is open. To see this, for any $f\in E(a,b,k,j)$, we have
\[
\sup_{\theta\in\partial\mathsf{Q}(0,1)}\min_{1 \leqslant s\leqslant a}\{\sup_{z\in\overline{\mathsf{Q}}(s\cdot\theta,k)}d_Y(f(z),\mathsf{T}_{-s\cdot \theta}\psi_j(z))\}
<
{1}/{b}.
\]
Hence we can take a positive $\delta$ with
$$\delta<{1}/{b}-\sup_{\theta\in\partial\mathsf{Q}(0,1)}\min_{1 \leqslant s\leqslant a}\{\sup_{z\in\overline{\mathsf{Q}}(s\cdot\theta,k)}d_Y(f(z),\mathsf{T}_{-s\cdot \theta}\psi_j(z))\},$$ 
and obtain an open neighbourhood of $f$
\[
N_{f,\delta,\overline{\mathsf{Q}}(0,5a\cdot k)}:=\{h\in\mathsf{Hol}(\mathbb{C}^n,Y):\sup_{z\in\overline{\mathsf{Q}}(0,5a\cdot k)} d_Y(h(z),f(z))<\delta\},
\]
which is clearly a subset of $E(a,b,k,j)$. 

Next, we claim that each $\cup_{a\geqslant 1}E(a,b,k,j)$ is dense. We only need to show that,
for any given holomorphic map $g: \mathbb{C}^n\rightarrow Y$ and $R>0$, there exists
some $F\in \cup_{a\geqslant 1}E(a,b,k,j)$ such that 
\begin{equation}\label{in neighborhood of g R}
\sup_{z\in \overline{\mathsf{Q}}(0,R)} d_Y(F(z),g(z))<{1}/{R}, \end{equation}

By the uniform continuity  of $\psi_j$ on compact sets, we can find a sufficiently small positive $\delta\ll 1$ such that
\begin{equation}\label{unifrom conti}
   d_Y(\psi_j(z_1),\psi_j(z_2))<{1}/{2b}\qquad{\scriptstyle(\forall\,z_1\,,\,z_2\,\in\,\overline{\mathsf{Q}}(0,k+\delta),\,\text{with}\,z_2\,\in\,\overline{\mathsf{Q}}(z_1\,,\,\delta)\,)}.
\end{equation}

Applying Lemma \ref{placing cubes} for the $\delta$ and $k\in\mathbb{N}$, we obtain a sufficiently large number $r>\lceil R \rceil+k$ satisfying that
there are finitely many pairwise disjoint closed hypercubes $\{\overline{\mathsf{Q}}(a_i,k+\delta)\}_{i\in\Gamma_{r,k,\delta}}$ out side $\mathsf{Q}(0,r)$ such that for any $\theta\in\partial\mathsf{Q}(0,1)$, there is some $m\in\mathbb{N}$ and $l\in\Gamma_{r,k,\delta}$ 
with
\begin{equation}\label{small deviation2}
    \mathsf{Q}(m\cdot\theta,k)\subset\overline{\mathsf{Q}}( a_l,k+\delta).
\end{equation}
On some small neighborhood of the union of the disjoint hypercubes
\[
\cup_{i\in\Gamma_{r,k,\delta}}{\overline{\mathsf{Q}}(a_i,k+\delta)}\cup\overline{\mathsf{Q}}(0,R),
\]
which is polynomially convex by Lemma~\ref{placing cubes},
we define
\begin{equation}\label{def of G}
    G(z)=\begin{cases}
    g(z), \qquad &\text{for}\,\,z\,\,\text{near}\,\,\overline{\mathsf{Q}}(0,R),\\
    \psi_j(z-a_i),\qquad&\text{for}\,\,z\,\,\text{near}\,\,\overline{\mathsf{Q}}(a_i,k+\delta),\,\forall i\,\in\Gamma_{r,k,\delta}.
\end{cases}
\end{equation}
 By the BOPA property of the Oka manifold $Y$ (\cite[p.~258]{Forstneric-Oka-book} and Remark~\ref{why can bopa}), we can find some holomorphic map $F: \mathbb{C}^n\rightarrow Y$  approximating $G$ within tiny deviation
\begin{equation}\label{approx G}
d_Y(F(z), G(z))
<
\min\{{1}/{b}, {1}/{R}\}/2
\qquad
{\scriptstyle
	(\,\forall\, z\,\in\, \bigcup_{i\in\Gamma_{r,k,\delta}}\, {\overline{\mathsf{Q}}(r_d\cdot \theta_{d,l},k+\delta)}\,\cup\,\overline{\mathsf{Q}}(0,R))}.
\end{equation}
Then \eqref{in neighborhood of g R} is automatically satisfied, and $F\in \cup_{a\geqslant 1}E(a,b,k,j)$,  because
\begin{align*}
    \sup_{z\in\overline{\mathsf{Q}}(0,k)}d_Y(\mathsf{T}_{m\cdot\theta}F(z),\psi_j(z))
    =&
    \sup_{z\in\overline{\mathsf{Q}}(m\cdot\theta,k)}d_Y(F(z),\psi_j(z-m\cdot\theta))\\
    [\text{see~\eqref{small deviation2}}]\qquad
    \leqslant&
    \sup_{z\in\overline{\mathsf{Q}}(a_l,k+\delta)}d_Y(F(z),\psi_j(z-a_l))\\
    &\qquad+\sup_{z\in\overline{\mathsf{Q}}(a_l,k+\delta)}d_Y(\psi_j(z-m\cdot\theta),\psi_j(z-a_l))\\
    [\text{use~\eqref{unifrom conti}}]\qquad
    \leqslant&
    \sup_{z\in\overline{\mathsf{Q}}(a_l,k+\delta)}d_Y(F(z),\psi_j(z-a_l))+{1}/{2b}\\
    [\text{check~\eqref{def of G}}]\qquad
    =&
    \sup_{z\in\overline{\mathsf{Q}}(a_l,k+\delta)}d_Y(F(z),G(z))+{1}/{2b}\\
    [\text{by~\eqref{approx G}}]\qquad
    <& \quad {1}/{b}.
\end{align*}
\qed



\section{\bf A harmonic analogue of Theorem~A}\label{Proof of theorem D}

A smooth function $f\in \mathscr{C}^{\infty}(\mathbb{R}^n)$
of $n$ variables $x_1, \dots, x_n$ is called harmonic, if it is in the kernel of the Laplacian operator
$
\Delta_x=\frac{\partial^2}{\partial x^2_1}+\cdots+\frac{\partial^2}{\partial x^2_n}
$. A classical result states that every harmonic function is real analytic. 
Denote by
\begin{equation}
    \label{define H(Rn, Rm)}
    \mathcal{H}(\mathbb{R}^n,\mathbb{R}^m)=\{f=(f_1, \dots, f_m):\mathbb{R}^n\rightarrow\mathbb{R}^m\,\,\,
:\,\,\,\Delta f_i=0,\, 1\leqslant i \leqslant m\}
\end{equation} 
the space of 
harmonic  maps
between two Euclidean spaces $\mathbb{R}^n$ and $\mathbb{R}^m$,
 equipped with the compact-open topology.
 Here is an analogue of Theorem~A, which generalizes a  result of~\cite{gomez2010slow}. 

\begin{thm}
\label{harmonic version}
For arbitrary 
 countable  elements $\{e^{[i]}\}_{i\geqslant 1}$ in the unit sphere $\mathbb{S}^{n-1}$ in $\mathbb{R}^n$, 
 the set 
 	\[
	\cap_{k\in\mathbb{N}}\cap_{ r\in\mathbb{R}_+}\mathsf{HC}(\mathsf{T}_{r\cdot e^{[k]}})\cap \mathsf{S}_{\phi}
	\] of 
	common hypercyclic harmonic maps  with slow growth 
 contains some image
 $\iota(G\cap \mathbb{B}(0,1))$,
 where
 $\iota: \mathbf{B}\hookrightarrow \mathcal{H}(\mathbb{R}^n,\mathbb{R}^m)$ is a continuous linear injective map from some Banach space $\mathbf{B}$, and where 
 $\mathbb{B}(0,1)\subset \mathbf{B}$ is the unit ball centered at $0$, and where
 $G\subset \mathbf{B}$ is some 
$G_\delta$-dense subset. 
\end{thm}

In the proof of Theorem~A, polynomials are used to approximate holomorphic functions. For  Theorem~\ref{harmonic version}, 
we will rely on harmonic polynomials $\mathsf{ker}\Delta_x\cap\mathbb{R}[x_1,\dots,x_n]$. 

By a construction in~\cite{gomez2010slow}, the $\mathbb{R}$-linear space of
harmonic polynomials
with $n\geqslant 2$ variables admits a  homogeneous basis $\{P_{i,j}(x)\}_{i,j\in\mathbb{N}}$  such that
\begin{equation}\label{harmonic basis}
\frac{\partial}{\partial x_1}\,P_{i+1,j}(x)\,
=\,
P_{i,j}(x)\qquad
{\scriptstyle
(\forall\,i,\,j\,\in\,\mathbb{N})
}.
\end{equation}
Thus we can define a linear isomorphism
	\begin{equation}\label{harmonic identification} 
	\varphi_x:
	\mathbb{R}^{\mathbb{N}^2}
	\longrightarrow
	\mathbb{R}[[x_1,\dots,x_n]]\cap \mathsf{ker} \Delta_x,
	\qquad
	(\alpha_{i,j})_{i,j\in\mathbb{N}}\mapsto
	\sum_{i,j\in\mathbb{N}} \alpha_{i,j} P_{i,j}(x),
	\end{equation}
making the following  diagram commutative
\[
\begin{tikzcd}
\mathbb{R}^{\mathbb{N}^2} \arrow[r, "{  B_1}"] \arrow[d, "\varphi_x"] & \mathbb{R}^{\mathbb{N}^2} \arrow[d, "\varphi_x"] \\
{\mathbb{R}[[x_1,\dots,x_n]]\cap \mathsf{ker}\Delta_x} \arrow[r, "  \frac{\partial}{\partial x_1} "]              & {\mathbb{R}[[x_1,\dots,x_n]]\cap \mathsf{ker}\Delta_x}.         
\end{tikzcd}
\]
	Here $B_1((\alpha_{i,j})_{i,j\in\mathbb{N}})=(\alpha_{i+1,j})_{i,j\in\mathbb{N}}$ is the backward shift operator with respect to the first index $i$.

Define $B_l=\varphi_x^{-1}
\circ\frac{\partial}{\partial x_l}\circ \varphi_x\,(2\leqslant l \leqslant n)$. For the same reason as~\eqref{diagram-1}, for any $(b_1, b_2, \dots, b_n)\in \mathbb{R}^n$, the following commutative diagram holds
\[
\begin{tikzcd}
\mathbb{R}^{\mathbb{N}^2} \arrow[rrr, "{\sum_{l=1}^n b_l B_l}"] \arrow[d, "\varphi_x"]                             &  &  & \mathbb{R}^{\mathbb{N}^2} \arrow[d, "\varphi_x"]     \\
{\mathbb{R}[[x_1,\dots,x_n]]\cap \mathsf{ker}\Delta_x} \arrow[rrr, "\sum_{l=1}^n b_l \frac{\partial}{\partial x_l}"] &  &  & {\mathbb{R}[[x_1,\dots,x_n]]\cap \mathsf{ker}\Delta_x}.
\end{tikzcd}
\]
We have a harmonic analogue of Lemma~\ref{suitable weight},
where $\ell^1(v)$ (see~\eqref{ell^1(v)}) is defined over $\mathbb{R}$ and  $I=\mathbb{N}^2$.

\begin{pro}\label{harmonic suitable weight}
    For any continuous function $\phi: \,\mathbb{R}_{\geqslant 0}\, \rightarrow\, \mathbb{R}_{> 0}$ of the shape~\eqref{growth phi}, 
    and homogeneous polynomials $\{P_{i,j}\}_{i,j\in \mathbb{N}}$ forming a basis of $\mathbb{R}[x_1,\dots,x_n]\cap \mathsf{ker}\Delta_x $, one can choose
    a positive weight $v=(v_{i,j})_{i,j\in \mathbb{N}}$ such that all the following considerations hold true simultaneously.
\begin{itemize}
\smallskip
\item[$\diamond \, 1.$]
For any $b=(b_1, b_2, \dots, b_n)\in\mathbb{R}^n\setminus\{\mathbf{0}\}$, the operator $T=
	\bigoplus_{u=1}^m\,e^{\sum^n_{l=1}b_l B_l}$ is continuous on $\big(\bigoplus_{u=1}^m\ell^1(v)\big)_{\ell^2}$ (recall notation~\eqref{infinite sum l^p}).
\smallskip
\item[$\diamond \,2.$]
 The $m$-fold direct sum of the map $\varphi_x$ 
\begin{align*}
    \bigoplus_{u=1}^m\varphi_x:
	\big(\bigoplus_{u=1}^m\ell^1(v)\big)_{\ell^2}
	&\rightarrow
	\mathcal{H}(\mathbb{R}^n,\mathbb{R}^m)\\
	\Big((\alpha^u_{s,t})_{s,t\in \mathbb{N}}\Big)_{1\leqslant u\leqslant m}
	&\mapsto
	\Big(\sum_{s,t\in \mathbb{N}} \alpha^u_{s,t}\, P_{s,t}(x)\Big)_{1\leqslant u\leqslant m}
\end{align*}
is well-defined, continuous,  with dense image. 
\smallskip
\item[$\diamond \, 3.$]
For any $\alpha\in\big(\bigoplus_{u=1}^m\ell^1(v)\big)_{\ell^2}$, one has
	\[
		\norm{ \bigoplus_{u=1}^m\varphi_x
		(\alpha)(x)}_{\mathbb{R}^m}
		\leqslant 
		\norm{\alpha}_{\big(\bigoplus_{u=1}^m\ell^1(v)\big)_{\ell^2}}
		\cdot
		\phi(\norm{x})
		\qquad
		{\scriptstyle
			(\forall\,x\,\in\, \mathbb{R}^n)
		}.
	\]
\end{itemize}
\end{pro}

In the holomorphic case,
the $\diamondsuit 1$ of Lemma~\ref{suitable weight} is guaranteed by the condition~\eqref{continuous condition}. Now, in the harmonic case, we will need a more delicate requirement~\eqref{harmonic v}
in the induction process.

\begin{proof}
 We can
demand $\sup_{i\in\mathbb{N}}\frac{v_{i,j}}{v_{i+1,j}}\leqslant M<+\infty$ to ensure that $B_1$ is continuous on $\ell^1(v)$. However, making other $B_l$ continuous is  subtle. 

Let us first partition the given basis $\{P_{i,j}\}_{i,j\in \mathbb{N}}$
according to degrees as $\bigsqcup_{k\in\mathbb{Z}_{\geqslant 0}}\mathcal{P}^k$, where $\mathcal{P}^k=\{P_{i,j}(x):\deg P_{i,j}(x)=k\}$.
Set
\[
\kappa : \mathbb{N}\times \mathbb{N}\rightarrow \mathbb{N},\qquad (i,j)\mapsto \deg P_{i,j}(x) +1.
\]
Let $\{e_{i, j}\}_{i, j\in \mathbb{N}}$ be the basis of $\mathbb{R}^{\mathbb{N}^2}$
so that $(\alpha_{i,j})_{i,j\in\mathbb{N}}$ in~\eqref{harmonic identification} reads as
$\sum_{i, j\in \mathbb{N}^2}\alpha_{i,j}\cdot e_{i, j}$.
Note that each derivative operator
$
\frac{\partial}{\partial x_l}$
maps the $\mathbb{R}$-linear subspace $
\mathsf{span}_{\mathbb{R}}\mathcal{P}^{k+1}$ to $\mathsf{span}_{\mathbb{R}}\mathcal{P}^{k}$. Hence we can 
write
\[
B_l(e_{i,j})=\sum_{s, t 
:\,
\kappa(s,t)=\kappa(i,j)-1} 
a^{l,i,j}_{s,t}\, e_{s,t}
\qquad
		{\scriptstyle
			(\forall\,1\,\leqslant\, l\,\leqslant\, n; \,\forall\, i,\, j\,\in \,\mathbb{N})
		}
\]
for some real constants $a^{l,i,j}_{s,t}$.

We now choose $(v_{i,j})$ inductively by $\kappa$. First, assign $v_{i,j}=1$ for $\kappa(i,j)=1$.
Assuming all $(v_{s,t})$ with $\kappa(s,t)\leqslant k$ have been determined ($k=1, 2, 3, \dots$), we choose positive $(v_{i,j})$ with $\kappa(i,j)= k + 1$ by the law
\begin{equation}\label{harmonic v}
v_{i,j}
\geqslant
\max_{1\leqslant l\leqslant n}
\sum_{s,t:\,\kappa(s,t)= k} |a^{l,i,j}_{s,t}|\cdot v_{s,t}.
\end{equation}
This guarantees~$\diamond 1$. Indeed,
\begin{align*}
    \norm{B_l\big(\sum_{i,j} \alpha_{i,j}\cdot e_{i,j}\big)}_{\ell^1(v)}\,
    =&\,
    \norm{\sum_{i,j}\alpha_{i,j}\sum_{s,t:\,\kappa(s,t)=\kappa(i,j)-1} a^{l,i,j}_{s,t} \cdot e_{s,t}}_{\ell^1(v)}
    \\
    \leqslant&\,
    \sum_{i,j}|\alpha_{i,j}|\sum_{s,t:\,\kappa(s,t)=\kappa(i,j)-1} |a^{l,i,j}_{s,t}|\cdot v_{s,t}\\
    \leqslant&\,
    \sum_{i,j}|\alpha_{i,j}|\cdot v_{i,j}
    \\
    \,=&\,
    \norm{\sum_{i,j} \alpha_{i,j}\cdot e_{i,j}}_{\ell^1(v)}.
\end{align*}
Therefore, \eqref{harmonic v} ensures
$\norm{B_l}_{\ell^1(v)\rightarrow\ell^1(v)}\leqslant 1$. 
Hence by the Taylor expansion we can see
$$
\norm{e^{\sum_{l=1}^n b_l B_l}}_{\ell^1(v)\rightarrow\ell^1(v)}
\leqslant
\prod_{l=1}^n\sum_{k=0}^{+\infty}\frac{{\norm{b_l B_l}^k_{\ell^1(v)\rightarrow\ell^1(v)}}}{k!}
\leqslant 
\prod_{l=1}^n\sum_{k=0}^{+\infty}\frac{{{|b_l|}^k}}{k!}
=
e^{\sum_{l=1}^n |b_l|}.
$$

In order to show $\diamond 2$ and $\diamond 3$,   
we take an enumeration $\Tilde{\iota}
:
\mathbb{N}
\overset{\sim}{\longrightarrow}
\mathbb{N}^2$, 
such that if $\kappa(i,j)\geqslant \kappa(i^\prime,j^\prime)$, then $\Tilde{\iota}^{-1}(i,j)\geqslant\Tilde{\iota}^{-1}(i^\prime,j^\prime)$. This can be achieved since each $t\in\mathbb{N}$, $\{(i,j):\kappa(i,j)=t\}\,$ is finite.
Set
\[
\hat v_k:= v_{\iota(k)},\qquad \hat P_{k}(x):=P_{\iota(k)}(x).
\]
Since $\phi$ grows faster than any polynomial, we can find a sequence of positive radii $\{R_k\}_{k\in\mathbb{N}}$ such that
\begin{equation}\label{estiamte harmonic 1-1}
	\phi(\norm{x})
	\geqslant
	2^k\cdot
	\Big|\hat P_k(x)\Big|
	\qquad
	{\scriptstyle
		(\forall\,\norm{x}\,\geqslant\, R_k).
	}
\end{equation}
Besides~\eqref{harmonic v}, we moreover require that the positive sequence $(\hat v_k)_{k\in \mathbb{N}}$ satisfies 
\begin{equation}\label{estiamte harmonic 1-2}
	\hat v_k\geqslant 
	\max\Big\{\frac{1}{\phi(0)}
	\sup_{\norm{x}\leqslant R_k} \Big|\hat P_k(x)\Big|,
	\,\frac{1}{2^k}
	\Big\}
		\qquad
	{\scriptstyle
		(\forall\,k\,\in\,\mathbb{N})}.
\end{equation}
Note that
\eqref{harmonic v} and \eqref{estiamte harmonic 1-2} can be achieved simultaneously along our induction process, thanks to the well-chosen $\tilde{\iota}$.

By~\eqref{estiamte harmonic 1-1} and \eqref{estiamte harmonic 1-2}, the following estimate follows by much the same argument as that of Lemma~\ref{suitable weight}
\[
		\norm{ \bigoplus_{u=1}^m\varphi_x
		(\alpha)(x)}_{\mathbb{R}^m}
		\leqslant 
		\norm{\alpha}_{\big(\bigoplus_{u=1}^m\ell^1(v)\big)_{\ell^2}}
		\cdot
		\phi(\norm{x})
		\qquad
		{\scriptstyle
			(\forall\,\alpha\,\in\,(\bigoplus_{u=1}^m\ell^1(v))_{\ell^2};\,\,\forall\,x\,\in\, \mathbb{R}^n)
		}.
\]
Thus $\diamond 3$ holds true. 
Consequently, $\diamond 2$ follows directly, since the $\mathbb{R}$-linear space spanned by the basis $\{P_{m,n}\}_{m,n\in \mathbb{N}}$ is dense in $\mathbb{R}[x_1,\dots,x_n]\cap \mathsf{ker}\Delta_x $.
\end{proof}

\begin{lem}\label{harmonic hypercyclic exp}
	For any positive weight $v=(v_{i,j})_{i,j\in \mathbb{N}}$ satisfying \eqref{harmonic v}, for any $b\in\mathbb{R}^n\setminus\{\mathbf{0}\}$, the operator $T=
	\bigoplus_{u=1}^m e^{\sum^n_{l=1}b_l B_l}$ is topologically transitive on $\big(\bigoplus_{u=1}^m\ell^1(v)\big)_{\ell^2}$.
\end{lem}
\begin{proof}
We mimic the proof of Lemma \ref{needed version}. The key point is to construct the following quasi-conjugate relation for some auxiliary positive weight $w$
	\[
	\begin{tikzcd}
		\ell^1(w) \arrow[rr, "e^{B_1}"] \arrow[d, "F"]           &  & \ell^1(w) \arrow[d, "F"] \\
		\ell^1(v) \arrow[rr, "e^{\sum_{l=1}^n b_l B_l}"] &  & \ell^1(v).               
	\end{tikzcd}
	\]	
	
For $b=(b_1, b_2, \dots, b_n)\in \mathbb{R}^n\setminus \{\mathbf{0}\}$ with Euclidean norm $\norm{b}_{\mathbb{R}^n}=r>0$, by some orthonormal transformation $\phi: \mathbb{R}^n\rightarrow \mathbb{R}^n$, we can find a new coordinate system $y=(y_1, y_2, \dots, y_n)$ of $\mathbb{R}^n$ such that $b$ now reads as
$\phi(r, 0, \dots, 0)$.
A key point is that $\phi$ preserves the Laplacian operator
$$
(\phi_*)^{-1} \Delta_x=\Delta_y
=
\frac{\partial^2}{\partial y^2_1}+\cdots+\frac{\partial^2}{\partial y^2_n},$$ 
since $\phi$ preserves the Euclidean metric on $\mathbb{R}^n$. 
Whence we obtain
$$\phi^*:\mathbb{R}[[x_1,\dots,x_n]]\cap \mathsf{ker}\Delta_x\rightarrow\mathbb{R}[[y_1,\dots,y_n]]\cap \mathsf{ker}\Delta_y,\qquad p(x)\mapsto p\circ \phi(y).$$ 
Now we arrive at the following commutative diagram
\[
\begin{tikzcd}
\mathbb{R}^{\mathbb{N}^2} \arrow[r, "\varphi_x"] \arrow[d, "\sum_{l=1}^n b_l B_l"] & {\mathbb{R}[[x_1,\dots,x_n]]\cap \mathsf{ker}\Delta_x} \arrow[r, "\phi^*"] \arrow[d, "\sum_{l=1}^n b_l \frac{\partial}{\partial x_l}"] & {\mathbb{R}[[y_1,\dots,y_n]]\cap \mathsf{ker}\Delta_y} \arrow[d, "r\frac{\partial}{\partial y_1}"] & \mathbb{R}^{\mathbb{N}^2} \arrow[l, "\varphi_y"] \arrow[d, "r B_1"] \\
\mathbb{R}^{\mathbb{N}^2} \arrow[r, "\varphi_x"]                                   & {\mathbb{R}[[x_1,\dots,x_n]]\cap \mathsf{ker}\Delta_x} \arrow[r, "\phi^*"]                                                             & {\mathbb{R}[[y_1,\dots,y_n]]\cap \mathsf{ker}\Delta_y}                                            & \mathbb{R}^{\mathbb{N}^2} \arrow[l, "\varphi_y"].                 
\end{tikzcd}
\]
By Lemma \ref{construct quasiconjugate}, we can choose a positive weight $w=(w_{i,j})_{i,j\in\mathbb{N}}$ to make $\varphi_x^{-1}\circ \phi^{*-1}\circ\varphi_y:\ell^1(w)\rightarrow\ell^1(v)$ continuous.
Moreover we can demand that $(w_{i,j})_{i,j\in \mathbb{N}}$
is increasing with respect to the first index $i$ to make sure $e^{B_1}$ is a well-defined continuous linear operator on $\ell^1(w)$ (See the proof of Lemma~\ref{needed version}).

Since $\phi^*$ induces a bijection between $\mathbb{R}[x_1,\dots,x_n]\cap\mathsf{ker} \Delta_x$ and $\mathbb{R}[y_1,\dots,y_n]\cap\mathsf{ker} \Delta_y$, it is clear that the map
$\varphi_x^{-1}\circ \phi^{*-1}\circ\varphi_y$ has dense image.
By mimicking the proof of Lemma~\ref{needed version} using  Theorem~\ref{key thm of DSW97}, we can show that $e^{r B_1}$ is hypercyclic on $\ell^1(w)$. Hence by Proposition~\ref{hypercyclic comparison}, we conclude that $e^{\sum_{l=1}^n b_l B_l}$ is hypercyclic on $\ell^1(v)$.
\end{proof}

\begin{proof}[Proof of Theorem~D]
Choose a positive weight $v=(v_K)_{K\in I}$ as in Proposition~\ref{harmonic suitable weight}. Then for any countable $\{b^{[k]}=(b_l^{[k]})_{1\leqslant l \leqslant n}\}_{k=1}^{+\infty}\subset\mathbb{R}^n\setminus\{\mathbf{0}\}$, the operators $
\bigoplus_{u=1}^m e^{\sum^n_{l=1}b^{[k]}_l B_l}$ are continuous on $\big(\bigoplus_{u=1}^m\ell^1(v)\big)_{\ell^2}$.
It follows from \eqref{harmonic basis} and Proposition~\ref{harmonic suitable weight} that	
we have the following quasi-conjugate relation
\[
\begin{tikzcd}
\big(\bigoplus_{u=1}^m\ell^1(v)\big)_{\ell^2} \arrow[rrr, "\bigoplus_{u=1}^m e^{\sum^n_{l=1}b^{[k]}_l B_l}"] \arrow[d, "\bigoplus_{u=1}^m\varphi_z"] &  & & \big(\bigoplus_{u=1}^m\ell^1(v)\big)_{\ell^2} \arrow[d, "\bigoplus_{u=1}^m\varphi_z"] \\
\mathcal{H}(\mathbb{R}^n,\mathbb{R}^m) \arrow[rrr, "\mathsf{T}_{b^{[k]}}"]                                                                                            &  & & \mathcal{H}(\mathbb{R}^n,\mathbb{R}^m).                                                     
\end{tikzcd}
\] 

Lemma \ref{harmonic hypercyclic exp} guarantees that all $\bigoplus_{u=1}^m e^{\sum_{l=1}^n b^{[k]}_l B_l} : \big(\bigoplus_{u=1}^m\ell^1(v)\big)_{\ell^2}\rightarrow \big(\bigoplus_{u=1}^m\ell^1(v)\big)_{\ell^2}$ are hypercyclic. 
By Proposition \ref{countable direction},  $\cap_{k\in\mathbb{N}}\mathsf{HC}(\bigoplus_{u=1}^m e^{\sum_{l=1}^n b^{[k]}_l B_l})$ is a $G_\delta$-dense subset of $\big(\bigoplus_{u=1}^m\ell^1(v)\big)_{\ell^2}$.
Let $\mathsf{B}_1$ be the closed unit ball of $\big(\bigoplus_{u=1}^m\ell^1(v)\big)_{\ell^2}$. 

Then 
    by Proposition~\ref{harmonic suitable weight},
    for any $\alpha\in\cap_{k\in\mathbb{N}}\mathsf{HC}(\bigoplus_{u=1}^m e^{\sum_{l=1}^n b^{[k]}_l B_l})\cap \mathsf{B}_1$, we have
    \[
		\norm{\bigoplus_{u=1}^m\varphi_x(\alpha)(x)}_{\mathbb{R}^n}
		\leqslant 
		\phi(\norm{x})
		\qquad
		{\scriptstyle
			(\forall\,x\,\in\, \mathbb{R}^n)
		}.
	\]
 By Proposition \ref{hypercyclic comparison},
 \[
 (\bigoplus_{u=1}^m\varphi_x)
 \Big(
 \cap_{k\in\mathbb{N}}\mathsf{HC}\big(\bigoplus_{u=1}^m e^{\sum_{l=1}^n b^{[k]}_l B_l}\big)\cap \mathsf{B}_1
 \Big)
 \subset
 \cap_{k\in\mathbb{N}} \mathsf{HC}(\mathsf{T}_{b^{[k]}})\cap\mathsf{S}_\phi\subset\mathcal{H}(\mathbb{R}^n,\mathbb{R}^m).
 \]
 
 Lastly, we recall the following
 
 \medskip\noindent
 {\bf Fact}\,(\cite{conejero2007hypercyclic}){\bf .}
 {\it
   For any $\theta\in\mathbb{S}^{n-1}$ and $ b\in\mathbb{R}_+\cdot \theta$, one has
    $\mathsf{HC}(\mathsf{T}_{\theta})=\mathsf{HC}(\mathsf{T}_{b})\subset\mathcal{H}(\mathbb{R}^n,\mathbb{R}^m)$.
}

\medskip
Whence 
\[
\cap_{k\geqslant 1} \mathsf{HC}(\mathsf{T}_{b^{[k]}})\cap\mathsf{S}_\phi
\,=\,
	\cap_{k\geqslant 1}\cap_{ r\in\mathbb{R}_+} \mathsf{HC}(\mathsf{T}_{r\cdot b^{[k]}})\cap\mathsf{S}_\phi
\]
and we conclude the proof.
\end{proof}

\begin{rmk}
Let us mention that, by a Runge-type approximation theorem for harmonic functions (cf.~\cite[Theorem 1]{du1969runge}),
a harmonic analogue of  Theorem~\ref{all trans}~(\cite[Theorem 1]{amazing-theorem})
follows after~\cite[Theorem~1.1]{bayart2016common}, {\em i.e.},
the set of common hypercyclic elements
\[
\cap_{b\in\mathbb{R}^n\setminus\{0\}}\,
\mathsf{HC}(\mathsf{T}_b)
\,
\subset\,
\mathcal{H}(\mathbb{R}^n,\mathbb{R}^m)
\]
 contains a $G_\delta$-dense subset. 
 \end{rmk}
 
 \medskip
 
 Similarly, we have the following harmonic analogue of Theorem~D by a Runge-type approximation theorem (cf.~\cite[Theorem 1]{du1969runge}) and our delicate arrangement of hypercubes (see Sect.~\ref{Hypercubes arrangement}).

\begin{thm}
In $\mathcal{H}(\mathbb{R}^n,\mathbb{R}^m)$ with the compact-open topology, 
$\cap_{b\in\mathbb{R}^n\setminus\{\mathbf{0}\}}\,
\mathsf{FHC}(\mathsf{T}_b)$ is a dense subset. \qed
\end{thm} 

\begin{ques}
Given integers $n\geqslant 3$, $m\geqslant 1$, 
is there an analogue of Theorem~C for $\mathcal{H}(\mathbb{R}^n,\mathbb{R}^m)$? 
\end{ques}

\begin{center}
	\bibliographystyle{alpha}
	\bibliography{article}
\end{center}

\end{document}